\documentclass{article}
\usepackage{amsmath}
\usepackage{amssymb}
\usepackage{amsthm}
\usepackage{float}

\newtheorem{Property}{Property}
\usepackage{paralist}
\usepackage{graphicx}

\newcommand{\Includegraphics}[2]%
{\ifFigs\includegraphics[      width=#1]{#2}%
	\else  \includegraphics[draft,width=#1]{#2}\fi}

\begin{document}
\title{On the emergence and properties of weird quasiperiodic attractors }
\author{\textit{Laura Gardini}$^{1,2},$ 
	\textit{ Davide Radi}$^{2,3}$, 
	\textit{Noemi Schmitt}$^{4}$, \\ 
	\textit{Iryna Sushko}$^{3,5}$, 
	\textit{Frank Westerhoff}$^{4}$\\ 
	$^{1}${\small Dept of Economics, Society and Politics, University of Urbino
		Carlo Bo, Italy}\\ 
	$^{2}${\small Dept of Finance, V\v{S}B - Technical University of Ostrava,
		Ostrava, Czech Republic}\\ 
	$^{3}${\small Dept of Mathematics for Economic, Financial and Actuarial
		Sciences,} \\ {\small Catholic University of Milan, Italy} \\ 
	$^{4}${\small Dept of Economics, University of Bamberg, Germany}\\ 
	$^{5}${\small Inst. of Mathematics, NAS of Ukraine, Ukraine}}

\date{}
\maketitle

\textbf{Abstract}
\medskip

We recently described a specific type of attractors of two-dimensional
discontinuous piecewise linear maps, characterized by two discontinuity lines
dividing the phase plane into three partitions, related to economic
applications. To our knowledge, this type of attractor, which we call a weird
quasiperiodic attractor, has not yet been studied in detail. They have a
rather complex geometric structure and other interesting properties that are
worth understanding better. To this end, we consider a simpler map that can
also possess weird quasiperiodic attractors, namely, a 2D discontinuous
piecewise linear map $F$ with a single discontinuity line dividing the phase
plane into two partitions, where two different homogeneous linear maps are
defined. Map $F$ depends on four parameters -- the traces and determinants of
the two Jacobian matrices. In the parameter space of map $F$, we obtain
specific regions associated with the existence of weird quasiperiodic
attractors; describe some characteristic properties of these attractors; and
explain one of the possible mechanisms of their appearance. 
\medskip 

\medskip


\section{Introduction}

Two-dimensional (2D) \textit{piecewise linear continuous} maps generating
rotational dynamics have been extensively studied in recent decades. In
particular, many researchers have considered the family of conservative maps
defined as $M=$ $M_{L}:(x,y)\rightarrow(\tau_{L}x-y,x)$ for $x\leq0,$ and
$M=M_{R}:(x,y)\rightarrow(\tau_{R}x-y,x)$ for $x\geq0,$ depending on two real
parameters, $\tau_{R}$ and $\tau_{L},$ which are the traces of the two
Jacobian matrices (see \cite{Beardon}, \cite{Sivak}, \cite{Lagarias a},
\cite{Lagarias b}, \cite{Lagarias c}, \cite{Garcia}, \cite{Roberts}). Since
straight lines through the origin (also called rays) are mapped by map $M$
into straight lines trough the origin, the dynamics of map $M$ can be
described by a circle map with well-defined rotation number 
$\theta(\tau_{R},\tau_{L})$. It was shown that in the $(\tau_{R},\tau_{L})$-parameter
plane, there exist sausage-like resonance regions, corresponding to a rational
rotation number $\theta(\tau_{R},\tau_{L})$ of the associated circle map. For
parameters in those regions, each point of the $(x,y)$-phase plane is
periodic. When $\theta(\tau_{R},\tau_{L})$ is irrational, there is a
topological conjugacy between the circle map and a rigid rotation, leading to
quasiperiodic trajectories, that is, each point of the $(x,y)$-phase plane
belongs to a quasiperiodic trajectory (see, e.g., \cite{Sivak}, \cite{Garcia}).

A more general \textit{linearly homogeneous} map satisfying 
$g(\alpha X)=\alpha g(X)$ for all $X\in \mathbb{R}^{d}$ 
and $\alpha\geq0$ is considered in \cite{Simpson-20-b}. In the 2D case,
this class of maps is associated with the well-known border collision form
$G_{\mu}$ for $\mu=0.$ We denote this map as $G_{0}$. Recall that the 2D
border collision normal form is defined as
\begin{equation}
G_{\mu}=\left\{
\begin{tabular}
[c]{l}%
$G_{L}:(x,y)\rightarrow(\tau_{L}x+y+\mu,-\delta_{L}x)$ for $x\leq0$\\
$G_{R}:(x,y)\rightarrow(\tau_{R}x+y+\mu,-\delta_{R}x)$ for $x\geq0$%
\end{tabular}
\ \ \right.  \label{G}%
\end{equation}
This map depends on parameters $\tau_{L},$ $\delta_{L},$ $\tau_{R},$ and
$\delta_{R},$ which are the trace and determinant of the Jacobian matrix in
partition $L$ (with $x\leq0$) and partition $R$ (with $x\geq0$), respectively;
parameter $\mu\neq0$ is usually scaled out by fixing it as $\mu=1.$ Since map
$G_{\mu}$ is quite important not only for various applications, but also for
the bifurcation theory of nonsmooth maps, its dynamics has been studied by
many researchers (see, e.g., \cite{Nusse-92}, \cite{SG-08}, \cite{diBernardo},
\cite{Simpson Meiss}, \cite{Simpson}, \cite{Zhusu}). In \cite{Simpson-20-b},
map $G_{0}$ is considered, whose dynamics can be characterized by the values
on the unit circle $S^{1}=\{x\in\mathbb{R}^{2}:||x||=1\}$. 
Introducing the polar coordinates $(x,y)=(\rho\cos
(\theta),\rho\sin(\theta)),$ the action of map $G_{0}$ can be split into a
radial expansion described by a scalar function $\rho$ and a rotation on
circle $S^{1}.$ For $\delta_{R}\delta_{L}\leq0,$ the map is noninvertible, and
one of the aims of \cite{Simpson-20-b} is to show that even in this
low-dimensional case, it is not a trivial task to determine the stability of
the fixed point in the origin. For $\delta_{R}\delta_{L}>0,$ the map is
invertible, and the stability of the fixed point is not guaranteed even if
this fixed point is attracting for both linear maps $G_{L}$ and $G_{R}$. As a
result, there are regions in the parameter space associated with divergence
and so-called \textit{dangerous bifurcation} (see, e.g., \cite{Do 06},
\cite{Do 07}, \cite{Ganguli}, \cite{GAS 09}, \cite{AZSBSG}).

The dynamics of 2D piecewise linear \textit{homogeneous}
\textit{discontinuous} maps, which can be quite intricate, is little studied
as yet. Such maps appear in some applications (see, e.g., \cite{Kollar},
\cite{Gardini et al 25}, \cite{Gardini et al 25c}). In particular, in
\cite{Gardini et al 25} an economic model is considered where the phase space
is separated into three partitions by two vertical discontinuity lines. One
linear homogeneous map acts in the middle partition, while the other acts in
the two external partitions. The origin is a fixed point for both maps, but
since it belongs to the middle partition, it is a virtual fixed point for the
external linear map. In the cited work, the existence of attractors that have
a special \textit{weird} structure -- not observed in other maps -- is
highlighted. In \cite{Gardini et al 25} and \cite{Gardini et al 25c} these
attractors are called \textit{weird quasiperiodic attractors}.

The goal of our paper is twofold. First, we analyze a mechanism for the
appearance of weird quasiperiodic attractors. Second, we describe their basic
properties. To simplify the analysis, one could consider the 2D border
collision normal form $G_{0}$. However, in order to introduce discontinuity,
some offset must be added to $G_{0}$ that destroys the property of
homogeneity. Instead, we consider the 2D piecewise linear map $F$ defined by
the same linear functions as in map $G_{0},$ but with the border line shifted
from $x=0$ to $x=h,$ $h\neq0$. In this way, map $F$ remains homogeneous but
becomes discontinuous:
\begin{equation}
F=\left\{
\begin{tabular}
[c]{l}%
$F_{L}:(x,y)\rightarrow(\tau_{L}x+y,-\delta_{L}x)$ for $x<h$\\
$F_{R}:(x,y)\rightarrow(\tau_{R}x+y,-\delta_{R}x)$ for $x>h$%
\end{tabular}
\ \ \ \ \ \right.  \label{F}%
\end{equation}

\begin{figure}
\includegraphics[width=1\textwidth]
{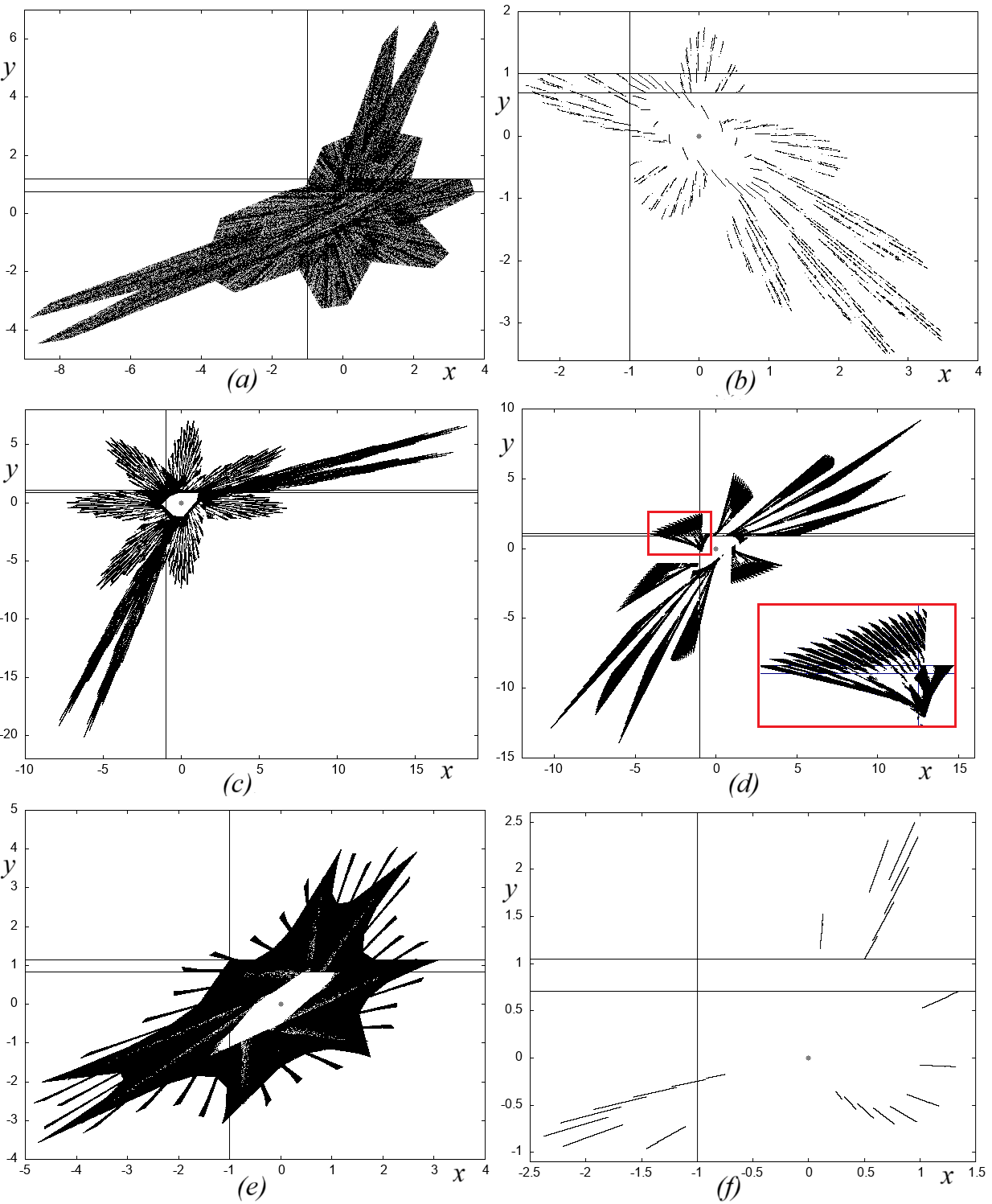} 
\caption{\label{f1} Attractors of map $F$ for (a) $\delta_{L}=0.75$, $\delta_{R}=1.2$,
$\tau_{L}=-0.7$, $\tau_{R}=-2.5$; (b) $\delta_{L}=0.7$, $\delta_{R}=1.001$,
$\tau_{L}=0.3$, $\tau_{R}=0.71$; (c) $\delta_{L}=0.9$, $\delta_{R}=1.1$,
$\tau_{L}=-2.5$, $\tau_{R}=-0.7$; (d) $\delta_{L}=0.9$, $\delta_{R}=1.1$,
$\tau_{L}=-2.5$, $\tau_{R}=-1.2$; (e) $\delta_{L}=0.84$, $\delta_{R}=1.15$,
$\tau_{L}=-1$, $\tau_{R}=-1.9$; (f) $\delta_{L}=1.05$, $\delta_{R}=0.7$,
$\tau_{L}=-0.75$, $\tau_{R}=-1.6$. Discontinuity line $x=-1$ and its images,
$y=\delta_{L}$ and $y=\delta_{R},$ by maps $F_{L}$ and $F_{R},$ respectively,
are also shown.}
\end{figure}

\begin{figure} 
\includegraphics[width=1\textwidth]%
{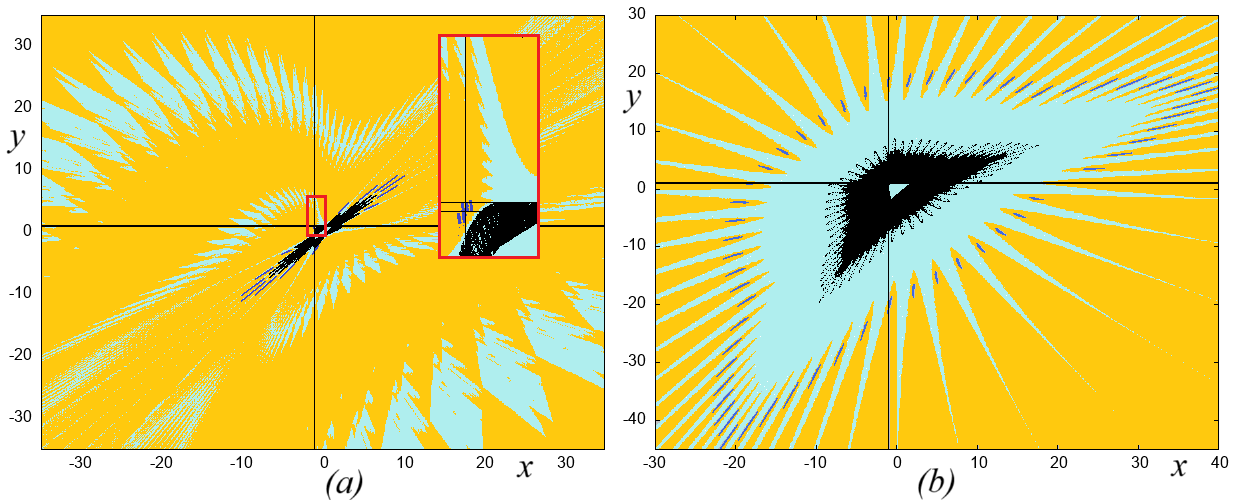}%
\caption{\label{f2} Two coexisting weird quasiperiodic attractors (shown in black and
dark blue) of map $F$ and their basins (in light blue and yellow) for (a)
$\delta_{L}=0.9$, $\delta_{R}=1.1$, $\tau_{L}=0.3$, $\tau_{R}=0.71;$ (b)
$\delta_{L}=0.9$, $\delta_{R}=1.11$, $\tau_{L}=-2$, $\tau_{R}=-1.91$. }
\end{figure}

Examples of the attractors of map $F$ for $h=-1,$ shown in Fig.~\ref{f1},
justify the use of the word \textit{weird}. Moreover, weird quasiperiodic
attractors can coexist, and in such cases their basins may also appear weird,
as shown in Fig.~\ref{f2}. At first glance, it is hard to believe that these
attractors are not chaotic. However, the characteristic property of map $F$ is
that both its components, maps $F_{L}$ and $F_{R},$ are homogeneous, with the
same fixed point in the origin, and it is quite straightforward to prove that
map $F$ cannot have hyperbolic cycles. This property facilitates the analysis
of its dynamics; in particular, one can immediately exclude the possibility
that its attractors are chaotic.\footnote{Among several known definitions of
chaos, we use the following one: a 2D map $F:I\rightarrow I,$ $I\subseteq \mathbb{R}^{2},$ 
is said to be chaotic on a closed invariant set $X\subseteq I$ if (a)
periodic points of $F$ are dense in $X$, and (b) $F$ is topologically
transitive, i.e., there is a dense aperiodic orbit on $X.$} Obviously, a
similar property holds also for a 2D piecewise linear map defined by several
nonhomogeneous linear maps having the same (nonzero) fixed point, which, via a
change of variables, can be translated to the origin. In our companion paper
\cite{Gardini et al 25b}, we discuss other classes of piecewise linear
discontinuous maps with weird quasiperiodic attractors. We think that similar
attractors may also be observed in other discontinuous maps, not only
piecewise linear ones.

In the present paper, we consider map $F$ given in (\ref{F}) as the simplest
representative of a class of 2D discontinuous piecewise linear maps, defined
by a finite number of linear homogeneous functions in different partitions of
the phase plane, separated by smooth discontinuity curves.

A weird quasiperiodic attractor (\textit{WQA} for short) $\mathcal{A}$ of map
$F$ is an attractor (a closed invariant attracting set with a dense
trajectory) with a rather complicated, say, weird, structure, which does not
include periodic points. Thus, $\mathcal{A}$ is neither an attracting cycle
nor a chaotic attractor: it is the closure of quasiperiodic trajectories. To
clarify, it is worth adding a few comments:

$\bullet$ By \textit{invariance} we mean that $F(\mathcal{A})=\mathcal{A}$.

$\bullet$ We use the standard (topological) definition of an attracting set,
namely, a closed invariant set $A$ is \textit{attracting} if there exists a
neighborhood $U(A)$ of $A$ such that 
$T(U(A))\subset U(A)$ and 
${\displaystyle\bigcap\limits_{i=0}^{\infty}}T^{i}(U(A))=A$. 
An attractor $\mathcal{A}$ is an
attracting set with a dense trajectory.

$\bullet$ As already mentioned, for map $F$ given in (\ref{F}), it is easy to
show that it cannot have hyperbolic cycles (see Property 5 in the next
section). This property is straightforward also for other piecewise linear
homogeneous maps.

$\bullet$ It is interesting to compare a WQA with other known kinds of
attractors that are also associated with quasiperiodic dynamics such as a
\textit{Cantor set attractor} (it can be observed, e.g., in a Lorenz map when
it is a gap map, see \cite{AGST2019}), or a \textit{critical attractor} (as,
e.g., in the logistic map at the Feigenbaum accumulation points, see
\cite{Sharkovsky}), or a \textit{strange nonchaotic attractor} (SNA for short,
observed in maps with quasiperiodic forcing, as described, e.g., in
\cite{Grebogi 84}, \cite{Feudel} for smooth maps or in \cite{Li 19}, \cite{Li
20}, \cite{Duan} for piecewise smooth maps). All these attractors have a
fractal structure. As for the structure of a WQA, the examples shown in Figs.
\ref{f1} and \ref{f2} have no fractal structure (see also other examples of
WQAs in the next sections). We think that it is a characteristic property of
the attractors of the considered class of maps, although it is not an easy
task to prove it in a generic case.

$\bullet$ In some special cases, map $F$ can have an attractor with a rather
simple structure (e.g., consisting of a finite number of segments), which can
be studied by means of a 1D map (a first return map on one of the segments).
We do not call such an attractor a WQA. For example, such a nongeneric case
may be observed for parameter values $\delta_{L}=0$ or $\delta_{R}=0,$ when
the entire partition defined by $x<0$ or $x>0,$ respectively, is mapped into
the straight line $y=0.$ As shown in our companion paper \cite{Gardini et al
25b}, if a 1D first return map to a suitable segment of this straight line
exists, with a finite number of discontinuity points, then this map is
topologically conjugated to a piecewise linear circle map with a rational or
irrational rotation number.

The paper is organized as follows. In the next section, several properties of
map $F$ are stated that are useful for describing its dynamics. In Sec. 3, we
present two examples of the bifurcation structure of the parameter space of
map $F,$ in the $(\tau_{L},\tau_{R})$-parameter plane for $\delta_{L}=0.9,$
$\delta_{R}=0.7,$ and in the $(\delta_{R},\tau_{R})$-parameter plane for
$\delta_{L}=0.9$, $\tau_{L}=-2.5$, with the regions related to the attracting
fixed point in the origin, weird quasiperiodic attractors, and divergence. We
obtain analytical expressions for the boundaries of the divergence regions in
the simplest case related to the basic and complementary to basic symbolic
sequences. In Sec. 4, we consider in detail the bifurcation structure near two
specific divergence regions and discuss, by means of several examples of the
phase portrait of map $F$, a mechanism of appearance of WQAs for parameters
values taken near those divergence regions. Some concluding remarks are
presented in Sec. 5.

\section{Preliminaries}

We consider a 2D discontinuous piecewise linear map 
$F:\mathbb{R}^{2}\rightarrow \mathbb{R}^{2},$ defined as follows:
\begin{equation}
F=\left\{
\begin{tabular}
[c]{ll}%
$F_{L}:\left(
\begin{array}
[c]{c}%
x\\
y
\end{array}
\right)  \rightarrow J_{L}\left(
\begin{array}
[c]{c}%
x\\
y
\end{array}
\right)  $ & if $x<-1$\\
$F_{R}:\left(
\begin{array}
[c]{c}%
x\\
y
\end{array}
\right)  \rightarrow J_{R}\left(
\begin{array}
[c]{c}%
x\\
y
\end{array}
\right)  $ & if $x>-1$%
\end{tabular}
\ \ \ \ \ \ \ \right.  \label{map}%
\end{equation}
where
\begin{equation}
J_{L}=\left(
\begin{array}
[c]{cc}%
\tau_{L} & 1\\
-\delta_{L} & 0
\end{array}
\right)  ,\ J_{R}=\left(
\begin{array}
[c]{cc}%
\tau_{R} & 1\\
-\delta_{R} & 0
\end{array}
\right)  \label{Jacobian}%
\end{equation}
Here, $\delta_{i}=\det J_{i},$ $\tau_{i}=trJ_{i},$ $i=L,R,$ are real
parameters. Note that we could choose any vertical line defined by $x=h,$
$h\neq0,$ as a discontinuity line, as parameter $h$ can be scaled out. Indexes
$L$ and $R$ in (\ref{map}) and (\ref{Jacobian}) refer to the left and right
partitions of the phase plane,
\[
D_{L}=\left\{  (x,y):x<-1\right\}  ,\ \ D_{R}=\left\{  (x,y):x>-1\right\}
\]
separated by the discontinuity line
\[
C_{-1}=\{(x,y):x=-1\}
\]
The images of $C_{-1}$ are denoted as
\[
C^{L}=F_{L}(C_{-1})=\{(x,y):y=\delta_{L}\},\ C^{R}=F_{R}(C_{-1}%
)=\{(x,y):y=\delta_{R}\}
\]
We call the discontinuity line $C_{-1},$ as well as its images and preimages,
\textit{critical lines} (of corresponding rank), following \cite{Mira-96}.

Let us state several properties of map $F$ that are easy to verify:

\begin{Property}
The \textit{fixed point} $O=(0,0)\in D_{R}$ of map $F_{R}$ is the unique fixed
point of map $F,$ provided no eigenvalue of $F_{R}$ equals $1,$ i.e.,
$1-\tau_{R}+\delta_{R}\neq0.$
\end{Property}

\begin{Property}
Invertibility of map $F$ is of

(a) $Z_{1}-Z_{0}-Z_{1}$ type for $0<\delta_{R}<\delta_{L},\ $or $\delta
_{L}<\delta_{R}<0;$

(b) $Z_{1}-Z_{2}-Z_{1}$ type for $\delta_{R}<\delta_{L}<0,$\ or $0<\delta
_{L}<\delta_{R};$

(c) $Z_{0}-Z_{1}-Z_{2}$ type for $\delta_{L}\delta_{R}<0;$

(d) $Z_{1}-Z_{\infty}-Z_{1}-Z_{0}$ type for $\delta_{L}=0,$ $\delta_{R}\neq0;$

(e) $Z_{0}-Z_{\infty}-Z_{0}-Z_{1}$ type for $\delta_{R}=0,$ $\delta_{L}\neq0;$

(f) $Z_{1}-Z_{2}-Z_{1}$ type for $\delta_{R}=\delta_{L}.$
\end{Property}

This property can be verified by considering the images of halfplanes $D_{L}$
and $D_{R},$ which define regions called \textit{zones} $Z_{j},$ where index
$j$ indicates the number of preimages of a point in zone $Z_{j}.$ In generic
cases (a), (b), or (c), zones $Z_{j}$ are separated by the critical lines
$C^{L}$ and $C^{R}$. In the nongeneric cases (d) and (e), $Z_{\infty}$ is the
critical line $C^{L}$ and $C^{R},$ respectively. Each point of these lines has
an infinite number of preimages (a straight line), since the complete
halfplane $D_{L}$ or $D_{R}$ is mapped into one straight line, $C^{L}$ or
$C^{R},$ respectively. In the nongeneric case (f), $Z_{2}$ is the critical
line $C^{L}=C^{R},$ whose points have two distinct preimages. For example, map
$F$ has $Z_{1}-Z_{2}-Z_{1}$ type of invertibility in Figs. \ref{f1}(a-e) and
$Z_{1}-Z_{0}-Z_{1}$ type of invertibility in Fig.~\ref{f1}(f).

\begin{Property}
In the $(\delta_{i},\tau_{i})$-parameter plane, $i=L,R,$ the eigenvalues
$\lambda_{1,2}^{i}=\frac{1}{2}(\tau_{i}\pm\sqrt{\tau_{i}^{2}-4\delta_{i}})$ of
matrix $J_{i}$ satisfy $\left\vert \lambda_{1,2}^{i}\right\vert <1$ in
the\textit{ }stability triangle $T^{i},$ defined as
\begin{equation}
T^{i}=\{(\delta_{i},\tau_{i}):\delta_{i}<1,\ -1-\delta_{i}<\tau_{i}%
<1+\delta_{i}\} \label{S}%
\end{equation}
This triangle is bounded by the segments of the straight lines $\tau
_{i}=1+\delta_{i}$ (related to $\lambda_{1}^{i}=1$), $\tau_{i}=-1-\delta_{i}%
$\ (related to $\lambda_{2}^{i}=-1$), and $\delta_{i}=1$ (related to
complex-conjugate eigenvalues $|\lambda_{1,2}^{i}|=1$).
\end{Property}

\begin{Property}
For the fixed point $O,$ the boundary of the stability domain $T^{R}$ defined
by $\tau_{R}=1+\delta_{R}$ corresponds to a degenerate $+1$ bifurcation, at
which any point of halfline\ $S^{R}=\{(x,y):x>-1,\ y=-\delta_{R}x\}$ is fixed,
while the boundary defined by $\tau_{R}=-1-\delta_{R}$ corresponds to a
degenerate flip bifurcation, at which any point (except the fixed point $O$)
of segment $S^{R^{2}}=\{(x,y):-1<x<1,\ y=\delta_{R}x\}$ is $2$-periodic. The
boundary defined by $\delta_{R}=1$ corresponds to a center bifurcation of the
fixed point $O,$ at which for a rational rotation number $\rho=m/n$ of matrix
$J_{R},$ i.e., for $\tau_{R}=2\cos(2\pi m/n),$ there exists an invariant
polygon with $n$ sides filled with nonhyperbolic cycles having rotation number
$m/n,$ while if $\rho$ is irrational, there then exists an invariant region
filled with quasiperiodic trajectories and bounded by an invariant ellipse
tangent to the discontinuity line $C_{-1}$.
\end{Property}

For more details on degenerate bifurcations in piecewise linear maps, we refer
to \cite{SG-10}.

Let $F_{\sigma}$ denote a composite map, $F_{\sigma}=F_{\sigma_{0}}\circ
F_{\sigma_{1}}...\circ F_{\sigma_{n-1}},$ where $\sigma=\sigma_{0}\sigma
_{1}...\sigma_{n-1}$ is a symbolic sequence with $\sigma_{j}\in\{L,R\},$
$n\geq2.$ The corresponding Jacobian matrix is denoted by $J_{\sigma}$, where
$J_{\sigma}={\displaystyle\prod\limits_{j=0}^{n-1}}J_{\sigma_{j}},$ 
and its characteristic polynomial is 
$P_{\sigma}(\lambda)=\lambda^{2}-trJ_{\sigma}\lambda+\det J_{\sigma}=0.$ It is easy to
see that the fixed point equation $F_{\sigma}(x,y)=(x,y)$ has a unique
solution $(x,y)=(0,0),$ provided that $P_{\sigma}(1)\neq0$ (i.e., matrix
$J_{\sigma}$ has no eigenvalue $1$), so that the unique hyperbolic fixed point
of $F_{\sigma}$ for any $\sigma$ is the fixed point $O.$ Thus, the following
property holds:

\begin{Property}
Map $F$ cannot have cycles of period $n$ for any $n\geq2,$ except for
nonhyperbolic $n$-cycles with an eigenvalue equal to $1$.
\end{Property}

Nonhyperbolic $n$-cycles with one eigenvalue equal to $1$ (not related to the
fixed point $O$), associated with condition $P_{\sigma}(1)=0,$ fill densely
specific segments, bounded or one-side unbounded. More precisely, let the
eigenvalues of matrix $J_{\sigma}$ be real and distinct (i.e., 
$(trJ_{\sigma})^{2}-4\det J_{\sigma}>0$), and let the corresponding eigenvectors be denoted
by $V^{\sigma}$ and $W^{\sigma}.$ Consider eigenvector $V^{\sigma}$ and its
$n-1$ images by map $F$ (all these sets obviously belong to straight lines
through the origin). Suppose the following \textit{admissibility} property is
satisfied: there exists a subset of $V^{\sigma}$ located in partition
$D_{\sigma_{0}}$ such that its $n-1$ images by $F$ are all located in the
proper partitions, according to the symbolic sequence $\sigma$. Let
$\{V_{j}^{\sigma}\}_{j=0}^{n-1}$ be a \textit{maximal admissible set}. This
means that there is a subset of $V^{\sigma}$ denoted by $V_{0}^{\sigma}$
belonging to $D_{\sigma_{0}}$ such that $F^{j}(V_{0}^{\sigma})=:$
$V_{j}^{\sigma}\subset D_{\sigma_{j}}$ for all $j=\overline{1,n-1},$ moreover,
one boundary point of (at least) one of the sets $V_{j}^{\sigma}$ belongs to
the discontinuity line $C_{-1}.$ The latter property implies that all the
other boundary points of the sets $V_{j}^{\sigma}$ are images of this point
(or these points). Now, let $P_{\sigma}(1)=0,$ i.e., one eigenvalue of matrix
$J_{\sigma}$ equals $1,$ and let $V^{\sigma}$ be the corresponding eigenvector
(for matrix $J_{\sigma}$ it is a straight line through the origin filled with
the fixed points). Then we can state the following

\begin{Property}
Map $F$ has a maximal admissible set of $n$ cyclic segments (bounded or
one-side unbounded), say, $S^{\sigma}=\{S_{j}^{\sigma}\}_{j=0}^{n-1},$
$n\geq2,$ filled with nonhyperbolic $n$-cycles (with eigenvalue $+1$) having
symbolic sequence $\sigma,$ provided that $P_{\sigma}(1)=0,$ 
$\det J_{\sigma}\neq1,$ and $S_{j}^{\sigma}=V_{j}^{\sigma}\subset D_{\sigma_{j}}$ for all
$j=\overline{0,n-1}.$
\end{Property}

Obviously, set $S^{\sigma}$ is not robust, i.e., it does not persist under
parameter perturbations. Put differently, the Lebesgue measure of the
parameter set satisfying $P_{\sigma}(1)=0$ is zero. Below, for short, we
denote a nonhyperbolic cycle having symbolic sequence $\sigma$ as 
$\sigma$-cycle. Examples of one-side unbounded segments 
$\{S_{j}^{\sigma}\}_{j=0}^{n-1}$ for various $\sigma$ can be seen in Fig.~\ref{D5Rexam}(b),(e),
Fig.~\ref{D5Rex3}(a), and Fig.~\ref{D5Lex1}(a), while bounded segments are
shown in Fig.~\ref{D5Lex2}(a) and Fig.~\ref{D5Lex5}(a).

Since map $F$ cannot have hyperbolic cycles, it cannot have chaotic attractors
either. Thus, we can state that in a generic case, a bounded attractor of map
$F$ is either the fixed point $O$ or a WQA. Note that set $S^{\sigma},$
consisting of $n$ cyclic segments filled with nonhyperbolic $\sigma$-cycles
(see Property 6), may be an attracting set, but it is not an attractor since
it does not include a dense trajectory. As mentioned in the Introduction, an
example of a nongeneric case may also be related to parameter values
$\delta_{R}=0$ or $\delta_{L}=0,$ when an attractor of map $F$ may be studied
by means of a 1D first return map (see \cite{Gardini et al 25b} for details).

With regard to the possible coexistence of attractors of map $F$, examples of
two coexisting WQAs are presented in Figs. \ref{f2} and \ref{D5Lex4}(b), while
examples of an attracting fixed point $O$ coexisting with a WQA are shown in
Figs. \ref{D5Rexam}(c),(f), \ref{D5Lex1}(b), and \ref{D5Lex2}(b).

As we will demonstrate in the next section, some boundaries of a region in the
parameter space, associated with WQAs of map $F$, are related to the existence
of specific sets $S^{\sigma}.$ As an example, we will consider \textit{basic
symbolic sequences}, $\sigma=LR^{n-1}$ and $\sigma=RL^{n-1},$ $n\geq3,$
assuming the simplest case of rotation number $\rho=1/n,$ as well as symbolic
sequences that are\textit{ complementary} to these sequences, namely,
$\sigma=L^{2}R^{n-2}$ and $\sigma=R^{2}L^{n-2},$ respectively. Here,
considering a (nonhyperbolic) cycle with a symbolic sequence $\sigma,$ by its
rotation number $\rho=m/n,$ we mean that along the cycle the trajectory makes
$m$ turns around the origin in $n$ iterations.

\section{Bifurcation structure of map $F$: \\ weird quasiperiodic attractors and
divergence}

\begin{figure}
[ht]
\begin{center}
\includegraphics[width=1\textwidth]{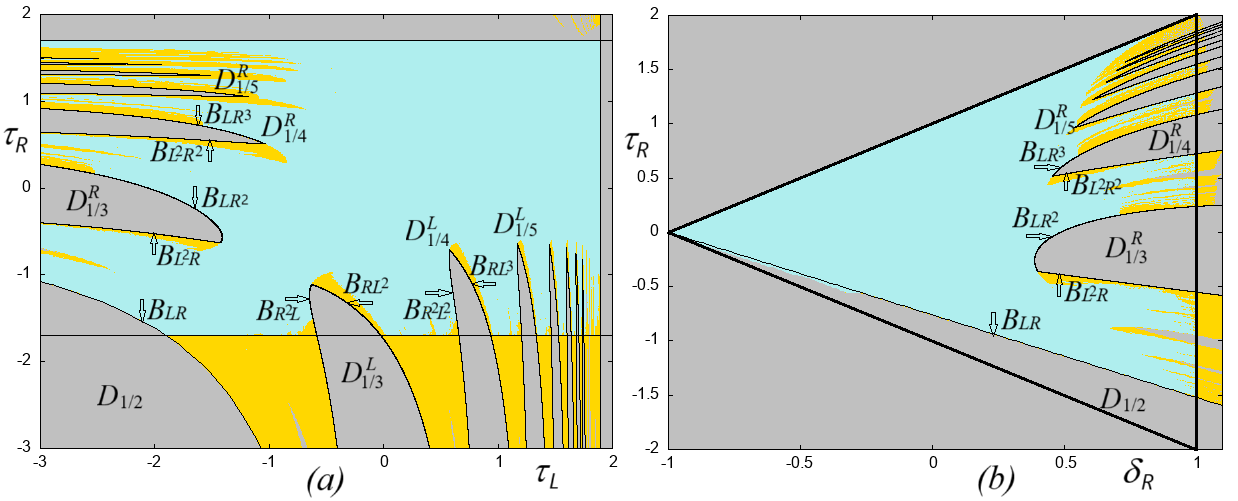}%
\caption{Bifurcation structure of map $F$ (a) in the 
	$(\tau_{L},\tau_{R})$-parameter plane for $\delta_{L}=0.9,$ $\delta_{R}=0.7$ , 
	and (b) in the $(\delta_{R},\tau_{R})$-parameter plane for $\delta_{L}=0.9,$ $\tau_{L}=-2.5$.
The boundaries $B_{LR^{n-1}}$ and $B_{L^{2}R^{n-2}}$ of the divergence region
$D_{1/n}^{R},$ defined in (\ref{B_LRn-1}) and (\ref{B_L2Rn-2}), respectively,
and the boundaries $B_{RL^{n-1}}$ and $B_{R^{2}L^{n-2}}$ of the divergence
region $D_{1/n}^{L}$ defined in (\ref{B_RLn-1}) and (\ref{B_R2Ln-2}),
respectively, are shown for $n=2,...,9.$ Blue and yellow parameter regions
indicate convergence to the fixed point $O$ and to a WQA, respectively. Gray
regions indicate divergence. }
\label{2D} 
\end{center}
\end{figure}

In Fig.~\ref{2D}(a), we present the bifurcation structure of map $F$ in the
$(\tau_{L},\tau_{R})$-parameter plane for $\delta_{L}=0.9,$ $\delta_{R}=0.7$.
Similarly, Fig.~\ref{2D}(b) shows the bifurcation structure of map $F$ in the
$(\delta_{R},\tau_{R})$-parameter plane for $\delta_{L}=0.9,$ $\tau_{L}=-2.5.$
In these figures, blue and yellow regions indicate convergence of an initial
point to the fixed point $O$ and to a WQA, respectively, while gray regions
indicate divergence. For parameter values as in Fig.~\ref{2D}(a), map $F$ has
$Z_{1}-Z_{0}-Z_{1}$ invertibility type since $0<\delta_{R}<\delta_{L}$ (see
Property 2(a)). In Fig.~\ref{2D}(b), map $F$ has $Z_{0}-Z_{1}-Z_{2}$
invertibility type for $\delta_{R}<0$ (Property 2(c)), $Z_{1}-Z_{0}-Z_{1}$
invertibility type for $0<\delta_{R}<0.9$ (Property 2(a)), $Z_{1}-Z_{2}-Z_{1}$
invertibility type for $\delta_{R}>0.9$ (Property 2(b)), and in special cases,
it is of $Z_{0}-Z_{\infty}-Z_{0}-Z_{1}$ invertibility type for $\delta_{R}=0$
(Property 2(e)), and $Z_{1}-Z_{2}-Z_{1}$ invertibility type for 
$\delta_{R}=\delta_{L}=0.9$ (Property 2(f)).

For $\delta_{R}=0.7,$ as in Fig.~\ref{2D}(a), the fixed point $O$ is
attracting in the horizontal strip $-1.7<\tau_{R}<1.7,$ while in Fig.
\ref{2D}(b), it is attracting for parameter values belonging to the indicated
stability triangle $T^{R}$ (see (\ref{S})). In both figures, the fixed point
$O$ can be globally attracting or coexisting either with a WQA, with divergent
trajectories, or with segments filled with nonhyperbolic cycles.

\begin{figure}
[ht]
\begin{center}
\includegraphics[width=1\textwidth]
{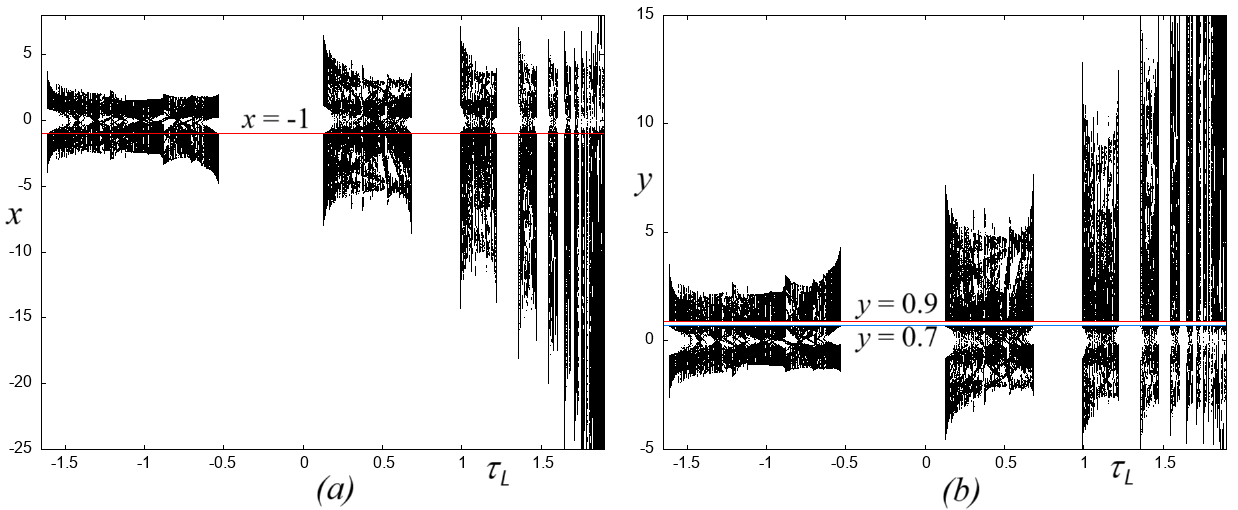}
\caption{A 1D bifurcation diagram of map $F$ showing $x$ versus $\tau_{L}$ in
(a), and $y$ versus $\tau_{L}$ in (b) for $\delta_{L}=0.9,$ $\delta_{R}=0.7,$
and $\tau_{R}=-2.$ }%
\label{1D}%
\end{center}
\end{figure}

In fact, the boundaries of the divergence regions that satisfy conditions
$P_{\sigma}(1)=0$ for specific $\sigma$ are related to the existence of $n$
cyclic \textit{halflines}, $S^{\sigma}=\{S_{j}^{\sigma}\}_{j=0}^{n-1},$
$n\geq2,$ filled with points of nonhyperbolic $\sigma$-cycles (see Property
6). Below we obtain the boundaries of the divergence regions, denoted by
$D_{1/n}^{R},$ associated with the basic and complementary to basic symbolic
sequences $\sigma=LR^{n-1}$ and $\sigma=L^{2}R^{n-2}$, as well as the
boundaries of the divergence regions, denoted by $D_{1/n}^{L},$ associated
with symbolic sequences $\sigma=RL^{n-1}$ and $\sigma=R^{2}L^{n-2}$. The lower
index $1/n$ means that in the considered case, the rotation number is
$\rho=1/n.$ The boundaries of regions $D_{1/n}^{R}$ and $D_{1/n}^{L}$ in Fig.
\ref{2D}(a) and the boundaries of regions $D_{1/n}^{R}$ in Fig.~\ref{2D}(b)
are shown for $n=2,...,9.$ In Fig.~\ref{1D}, the 1D bifurcation diagrams,
showing $x$ versus $\tau_{L}$ in (a) and $y$ versus $\tau_{L}$ in (b),
illustrate the dynamics of map $F$ for fixed $\tau_{R}=-2$ and increasing
$\tau_{L}.$ It can be seen that the divergence regions $D_{1/n}^{L}$ (with no
bounded dynamics) are crossed, and the regions between them correspond to
WQAs. Note that for $\tau_{R}=-2,$ the fixed point $O$ is a saddle. 

The divergence regions mentioned above are closely related to the so-called
\textit{dangerous bifurcations} occurring in the 2D border collision normal
form $G_{\mu}$. These bifurcations were studied quite intensively due to their
importance in the applied context (see, e.g., \cite{Hassouneh}, \cite{Ganguli}%
, \cite{Do 06}, \cite{Do 07}, \cite{GAS 09}, \cite{AZSBSG}). Recall that in
the 2D border collision normal form $G_{\mu}$ defined in (\ref{G}), a
dangerous bifurcation occurs at $\mu=0$ if for $\mu<0$ and for $\mu>0$ an
attractor or several attractors of map $G_{\mu}$ coexist with divergent
trajectories (i.e., with an attractor at infinity). In this case, as $\mu$
tends to $0,$ any attractor and its basin decrease in size linearly with
respect to $\mu$, and at $\mu=0,$ they shrink to point $(x,y)=(0,0).$ At
$\mu=0,$ the generic trajectory thus diverges. As a result, there exists one
or several attractors before and after the bifurcation, i.e., for $\mu<0$ and
$\mu>0,$ but near $\mu=0,$ the basins of bounded trajectories become quite
small, so that even small perturbations may lead to divergence, motivating the
name \textit{dangerous} bifurcation.

A necessary condition for a dangerous bifurcation to occur is the coexistence
of a bounded attractor with an attractor at infinity, which can be studied
using compactification of the phase plane (see, e.g., \cite{Simpson 16}).
Recall that in the parameter space of map $G_{\mu},$ the boundaries of the
divergence regions related to dangerous bifurcations satisfy conditions
$P_{\sigma}(1)=0$ for specific symbolic sequences $\sigma,$ as well as the
corresponding admissibility conditions. Crossing the boundary of a divergence
region defined by $P_{\sigma}(1)=0$ for some fixed $\mu\neq0,$ an $n$-cycle of
map $G_{\mu}$ with symbolic sequence $\sigma$ undergoes a \textit{degenerate
transcritical bifurcation}: its points tend to infinity while an eigenvalue
tends to $1.$ For $\mu=0,$ a kind of \textit{degenerate }$+1$\textit{
bifurcation} occurs: there are $n$ stable and $n$ unstable invariant halflines
issuing from the fixed point in the origin, and these halflines are filled
with points of nonhyperbolic $\sigma$-cycles. We refer to \cite{Ganguli},
\cite{GAS 09}, and \cite{AZSBSG} for more details related to the dangerous
bifurcation in map $G_{\mu},$ occurring at $\mu=0.$ In any case, condition
$P_{\sigma}(1)=0$ depends only on the parameters of matrices $J_{L}$ and
$J_{R}$, i.e., from $\tau_{i}$ and $\delta_{i},$ $i=L,R.$ These matrices are
defined in the same way for both maps $G_{\mu}$ and $F$, confirming an
intuitive idea about the similarity of the dynamics of maps $G_{\mu}$ and $F$
at infinity. However, map $F$ is discontinuous, and its border line is not
$x=0$ but $x=-1,$ which leads to some specific new properties mentioned below.

Let us obtain the boundaries of the divergence regions (satisfying conditions
$P_{\sigma}(1)=0$) in the parameter space of map $F$ by considering the
simplest case, namely, regions $D_{1/n}^{R}$ associated with symbolic
sequences $\sigma=LR^{n-1}$ and $\sigma=L^{2}R^{n-2},$ $n\geq3,$ and regions
$D_{1/n}^{L}$ associated with $\sigma=RL^{n-1}$ and $\sigma=R^{2}L^{n-2}.$

\subsection{$\sigma=LR^{n-1},$ $n\geq3$}

Consider first a composite map $F_{\sigma}$ with $\sigma=LR^{n-1},$ $n\geq3.$
Matrix $J_{LR^{n-1}}=J_{R}^{n-1}J_{L}$ can be defined as follows:
\[
J_{LR^{n-1}}=\left(
\begin{array}
[c]{cc}%
a_{n-1} & a_{n-2}\\
-\delta_{R}a_{n-2} & -\delta_{R}a_{n-3}%
\end{array}
\right)  \left(
\begin{array}
[c]{cc}%
\tau_{L} & 1\\
-\delta_{L} & 0
\end{array}
\right)  = \] 

\[ \left(
\begin{array}
[c]{cc}%
\tau_{L}a_{n-1}-\delta_{L}a_{n-2} & a_{n-1}\\
-\tau_{L}\delta_{R}a_{n-2}+\delta_{L}\delta_{R}a_{n-3} & -\delta_{R}a_{n-2}%
\end{array}
\right)
\]
where $a_{k}$ is the solution of the second-order homogeneous linear
difference equation
\[
a_{k}=\tau_{R}a_{k-1}-\delta_{R}a_{k-2},\ \ a_{0}=1,\ \ a_{1}=\tau
_{R},\ k=2,3,...
\]
The eigenvalues
\[
\lambda_{1,2}^{LR^{n-1}}=0.5(\tau_{L}a_{n-1}-(\delta_{L}+\delta_{R})a_{n-2}%
\pm\sqrt{(\tau_{L}a_{n-1}-(\delta_{L}+\delta_{R})a_{n-2})^{2}-4\delta
_{L}\delta_{R}^{n-1}})
\]
of matrix $J_{LR^{n-1}}$ are real if $(\tau_{L}a_{n-1}-(\delta_{L}+\delta
_{R})a_{n-2})^{2}-4\delta_{L}\delta_{R}^{n-1}>0,$ and the corresponding
eigenvectors have slopes
\begin{equation}
K_{i}^{LR^{n-1}}=\frac{\delta_{R}(\delta_{L}a_{n-3}-\tau_{L}a_{n-2})}%
{\delta_{R}a_{n-2}+\lambda_{i}},\ i=1,2 \label{K_LRn-1}%
\end{equation}
Transition from real to complex-conjugate eigenvalues occurs crossing a
parameter set denoted $E_{LR^{n-1}}$ and defined as follows:
\begin{equation}
E_{LR^{n-1}}:\ \ (\tau_{L}a_{n-1}-(\delta_{L}+\delta_{R})a_{n-2})^{2}%
-4\delta_{L}\delta_{R}^{n-1}=0 \label{E_LRn-1}%
\end{equation}

The characteristic polynomial of $J_{LR^{n-1}}$ can be written as follows:
\begin{equation}
P_{LR^{n-1}}(\lambda)=\lambda^{2}+\left(  -\tau_{L}a_{n-1}+(\delta_{L}%
+\delta_{R})a_{n-2}\right)  \lambda+\delta_{L}\delta_{R}^{n-1} \label{P_LRn-1}%
\end{equation}

Considering the linear map $F_{LR^{n-1}},$ i.e., map 
$F_{LR^{n-1}}:(x,y)\rightarrow J_{LR^{n-1}}(x,y),$ without its relation to the piecewise
linear map $F$, condition $P_{LR^{n-1}}(1)=0$ (for 
$\delta_{L}\delta_{R}^{n-1}\neq1$) means that map $F_{LR^{n-1}}$ has an eigenvalue equal to $1$ and
the straight line through the origin, say $y=K^{LR^{n-1}}x,$ associated with
the corresponding eigenvector, is filled with fixed points, where
\begin{equation}
K^{LR^{n-1}}=\frac{\delta_{R}(\delta_{L}a_{n-3}-\tau_{L}a_{n-2})}{\delta_{R}a_{n-2}+1} \label{K}%
\end{equation}
However, considering the piecewise linear map $F$, condition
\begin{equation}
P_{LR^{n-1}}(1)=1-\tau_{L}a_{n-1}+(\delta_{L}+\delta_{R})a_{n-2}+\delta_{L}\delta_{R}^{n-1}=0 \label{P_LRn-1(1)}%
\end{equation}
is necessary but not sufficient for the existence of set $S^{LR^{n-1}}$
consisting of $n$ segments filled with points of nonhyperbolic 
$LR^{n-1}$-cycles. According to Property 6, the admissibility conditions must be
satisfied, that is, these segments must be located in the proper partitions.
Here we are interested in the boundaries of the divergence regions, so these
segments are \textit{halflines} and their "end points at infinity" are, for
the compactified map, periodic as well, with symbolic sequence $LR^{n-1}.$ In
the considered case, denoting the starting halfline belonging to $D_{L}$ as
$S_{0}^{LR^{n-1}},$ so that
\[
S_{0}^{LR^{n-1}}=\{(x,y):y=K^{LR^{n-1}}x,\ x<-1\}
\]
(assuming $K^{LR^{n-1}}\neq\infty,$ that is, $\delta_{R}a_{n-2}+1\neq0$, see
(\ref{K}), that holds for $\delta_{R}\neq1$), its images 
$S_{j}^{LR^{n-1}}=F(S_{j-1}^{LR^{n-1}})$ for $j=\overline{1,n-1}$ must be completely located
in partition $D_{R}$. For short, we denote these admissibility conditions as
$(A^{LR^{n-1}})$, that is,
\begin{equation}
(A^{LR^{n-1}}):\ \ S_{0}^{LR^{n-1}}\subset D_{L},\ S_{j}^{LR^{n-1}}\subset
D_{R},\ \text{ }j=\overline{1,n-1} \label{(ALRn-1)}%
\end{equation}
(an example of admissible set $S^{LR^{4}}$ consisting of five halflines filled
with nonhyperbolic $LR^{4}$-cycles is shown in Fig.~\ref{D5Rexam}(a)).

In contrast to the map $G_{0}$, map $F$ can also have a set $S^{LR^{n-1}}$
satisfying condition $P_{LR^{n-1}}(1)=0$ and admissibility conditions
$(A^{LR^{n-1}})$, but consisting of \textit{bounded segments} (see Fig.
\ref{D5Rex3}(a) for an example), and in such a case the corresponding
parameter set is not associated with the boundary of a divergence region (see
the point marked $g$ in Fig.~\ref{D5R}(a) corresponding to Fig.~\ref{D5Rex3}%
(a)). In fact, comparing with map $G_{0}$, in map $F$ the border line is
shifted to the left, from $x=0$ to $x=-1,$ so that partition $D_{R}$ is
enlarged, providing additional space for the segments of set $S^{LR^{n-1}}$ to
be still admissible, which in the case of map $G_{0}$ would be not admissible
(or, in other words, virtual).

To summarize, we can state that among the parameter values satisfying
$P_{LR^{n-1}}(1)=0$ given in (\ref{P_LRn-1(1)}), there may be a subset denoted
$B_{LR^{n-1}}$ that satisfies the admissibility conditions $(A^{LR^{n-1}})$
given in (\ref{(ALRn-1)}):
\begin{equation}
B_{LR^{n-1}}:\ \ 1-\tau_{L}a_{n-1}+(\delta_{L}+\delta_{R})a_{n-2}+\delta
_{L}\delta_{R}^{n-1}=0\text{ and }(A^{LR^{n-1}})\text{ holds} \label{B_LRn-1}%
\end{equation}
where segments $S_{j}^{LR^{n-1}},\ $ $j=\overline{1,n-1},$ in $(A^{LR^{n-1}})$
can be either one-side unbounded (then $B_{LR^{n-1}}$ belongs to a boundary of
the divergence region $D_{1/n}^{R}$), or bounded.

\subsection{$\sigma=L^{2}R^{n-2},$ $n\geq3$}

Following analogous steps for the complementary symbolic sequence
$\sigma=L^{2}R^{n-2},$ we obtain a set of parameter values, denoted by
$B_{L^{2}R^{n-2}},$ satisfying $P_{L^{2}R^{n-2}}(1)=0,$ and the corresponding
admissibility conditions:
\[
B_{L^{2}R^{n-2}}: \ \  1+(\tau_{L}(\delta_{L}+\delta_{R})-\delta_{L}\tau
_{R})a_{n-3}-(\tau_{L}^{2}-2\delta_{L})a_{n-2}+\delta_{L}^{2}\delta_{R}%
^{n-2}=0 
\] 
\begin{equation} 
\text{ and }(A^{L^{2}R^{n-2}})\text{ holds} 
\label{B_L2Rn-2}%
\end{equation}
where
\[
(A^{L^{2}R^{n-2}}):\ \ S_{0}^{L^{2}R^{n-2}}\subset D_{L},\ S_{1}^{L^{2}%
R^{n-2}}\subset D_{L},\ S_{j}^{L^{2}R^{n-2}}\subset D_{R},\text{ }%
j=\overline{2,n-1}%
\]
and
\[
S_{0}^{L^{2}R^{n-2}}=\{(x,y):y=K^{L^{2}R^{n-2}}x,\ x<-1\}
\]
with
\[
K^{L^{2}R^{n-2}}=\frac{\delta_{R}(\delta_{L}a_{n-3}+\tau_{L}\left(  \delta
_{L}a_{n-4}-\tau_{L}a_{n-3}\right)  )}{\delta_{R}(\tau_{L}a_{n-3}-\delta
_{L}a_{n-4})+1}%
\]
(see an example of set $S^{L^{2}R^{3}}$ consisting of five halflines filled
with nonhyperbolic $L^{2}R^{3}$-cycles in Fig.~\ref{D5Rexam}(d)).

In this way, two boundaries, $B_{LR^{n-1}}$ and $B_{L^{2}R^{n-2}},$ of the
divergence region $D_{1/n}^{R},$ $n\geq3,$ are obtained (see Fig.~\ref{2D}).

Note that for a parameter point belonging to region $D_{1/n}^{R},$ map $F$ has
an attracting $n$-cycle at infinity with symbolic sequence either $LR^{n-1}$
or $L^{2}R^{n-2},$ associated with the unstable eigenvector of matrix
$J_{LR^{n-1}}$ or $J_{L^{2}R^{n-2}},$ respectively, and its $n-1$ images by
map $F$ (see, e.g., vectors $V_{j}^{LR^{4}}$ and $V_{j}^{L^{2}R^{3}},$
$j=\overline{0,4}$, in Fig.~\ref{D5Rexam}(a) and (d), respectively). A
boundary, say $H_{LR^{n-1}},$ separating the corresponding subregions of
region $D_{1/n}^{R}$ can be defined by the condition that the eigenvector
of\ matrix $J_{LR^{n-1}}$ has slope $0$ (and, thus, the slope of its preimage
by $F_{R}$ is equal to infinity). Then, for a parameter point belonging to
region $D_{1/n}^{R},$ we have the following:

$\bullet$ on one side of boundary $H_{LR^{n-1}},$ it holds that the unstable
eigenvectors $V_{j}^{LR^{n-1}}$ are admissible, while eigenvectors
$V_{j}^{L^{2}R^{n-2}}$ are not admissible (so that an attracting 
$LR^{n-1}$-cycle at infinity exists, while the $L^{2}R^{n-2}$-cycle is virtual);

$\bullet$ on the other side of boundary $H_{LR^{n-1}}$, the unstable
eigenvectors $V_{j}^{L^{2}R^{n-2}}$ are admissible, while eigenvectors
$V_{j}^{LR^{n-1}}$ are not admissible (so that an attracting 
$L^{2}R^{n-2}$-cycle at infinity exists, while the $LR^{n-1}$-cycle is virtual).

So, assuming that in (\ref{K_LRn-1}) $\delta_{R}\neq0,$ from 
$K_{i}^{LR^{n-1}}=0$ boundary $H_{LR^{n-1}}$ can be defined as follows:
\begin{equation}
H_{LR^{n-1}}:\ \ \delta_{L}a_{n-3}-\tau_{L}a_{n-2}=0,\ \ H_{LR^{n-1}}\subset D_{1/n}^{R} 
\label{HR}%
\end{equation}
Transition from real to complex-conjugate eigenvalues of matrix 
$J_{L^{2}R^{n-2}}$ (when map $F$ has no cycles at infinity) occurs crossing a parameter
set
\begin{equation}
E_{L^{2}R^{n-2}}:\ ((\delta_{L}\tau_{R}-\tau_{L}(\delta_{L}+\delta
_{R}))a_{n-3}+(\tau_{L}^{2}-2\delta_{L})a_{n-2})-4\delta_{L}^{2}\delta
_{R}^{n-2}=0 \label{E_L2Rn-2}%
\end{equation}
As an example, region $D_{1/5}^{R}$ is shown magnified in Fig.~\ref{D5R}(a),
where besides its boundaries $B_{LR^{4}}$ and $B_{L^{2}R^{3}},$ also curves
$H_{LR^{4}},$ $E_{LR^{4}},$ and $E_{L^{2}R^{3}}$ are shown, given in
(\ref{HR}), (\ref{E_LRn-1}), and (\ref{E_L2Rn-2}), respectively.

\subsection{$\sigma=RL^{n-1}$ and $\sigma=R^{2}L^{n-2},$ $n\geq3$}

Analogous boundaries associated with the divergence regions $D_{1/n}^{L}$ and
symbolic sequences $\sigma=RL^{n-1}$ and $\sigma=R^{2}L^{n-2}$ can be obtained
by changing $L\leftrightarrow R$ in all the above notations. In particular,
the boundary $B_{RL^{n-1}}$ of divergence regions $D_{1/n}^{L}$ can be defined
as
\begin{equation}
B_{RL^{n-1}}:\ \ 1-\tau_{R}b_{n-1}+(\delta_{R}+\delta_{L})b_{n-2}+\delta
_{R}\delta_{L}^{n-1}=0\text{ and }(A^{RL^{n-1}})\text{ holds} \label{B_RLn-1}%
\end{equation}
where the admissibility conditions $(A^{RL^{n-1}})$ are as follows:
\[
(A^{RL^{n-1}}):\ \ S_{0}^{RL^{n-1}}\subset D_{R},\ S_{j}^{RL^{n-1}}\subset
D_{L},\ \text{ }j=\overline{1,n-1}%
\]
and $b_{k}$ is a solution of the second-order homogeneous linear difference
equation%
\[
b_{k}=\tau_{L}b_{k-1}-\delta_{L}b_{k-2},\ \ b_{0}=1,\ \ b_{1}=\tau
_{L},\ \ k=2,3,...
\]
Transition from real to complex-conjugate eigenvalues of matrix $J_{RL^{n}}$
is associated with crossing a parameter set $E_{RL^{n-1}}$ satisfying
\begin{equation}
E_{RL^{n-1}}:\ \ (\tau_{R}b_{n-1}-(\delta_{L}+\delta_{R})b_{n-2})^{2}%
-4\delta_{R}\delta_{L}^{n-1}=0 \label{E_RLn}%
\end{equation}
The boundary $B_{R^{2}L^{n-2}}$ of divergence regions $D_{1/n}^{L}$ can be
defined as
\[
B_{R^{2}L^{n-2}}:\ \ 1+(\tau_{R}(\delta_{L}+\delta_{R})-\delta_{R}\tau
_{L})b_{n-3}-(\tau_{R}^{2}-2\delta_{R})b_{n-2}+\delta_{R}^{2}\delta_{L}%
^{n-2}=0 
\] 
\begin{equation} 
\text{ \ and }(A^{R^{2}L^{n-2}})\text{ holds} \label{B_R2Ln-2}%
\end{equation}
where
\[
(A^{R^{2}L^{n-2}}):\ \ S_{0}^{R^{2}L^{n-2}}\subset D_{R},\ S_{1}^{R^{2}%
L^{n-2}}\subset D_{R},\ S_{j}^{R^{2}L^{n-2}}\subset D_{L},\text{ }%
j=\overline{2,n-1}%
\]
Transition from real to complex-conjugate eigenvalues of matrix $J_{R^{2}%
L^{n-2}}$ are associated with a parameter set
\begin{equation}
E_{R^{2}L^{n-2}}:\ ((\delta_{R}\tau_{L}-\tau_{R}(\delta_{R}+\delta
_{L}))b_{n-3}+(\tau_{R}^{2}-2\delta_{R})b_{n-2})-4\delta_{R}^{2}\delta
_{L}^{n-2}=0 \label{E_R2Ln-2}%
\end{equation}
and boundary $H_{RL^{n-1}}$ can be defined as follows:
\begin{equation}
H_{RL^{n-1}}:\ \ \delta_{R}b_{n-3}-\tau_{R}b_{n-2}=0,\ \ H_{RL^{n-1}}\subset D_{1/n}^{L} 
\label{HL}%
\end{equation}
Figure \ref{D5L}(a) presents an example of the divergence region $D_{1/5}^{L}$
together with its boundaries $B_{RL^{4}}$ and $B_{R^{2}L^{3}},$ as well as the
curves $H_{RL^{4}},\ E_{RL^{4}},$ and $E_{R^{2}L^{3}}$, given in
(\ref{B_RLn-1}), (\ref{B_R2Ln-2}), (\ref{E_RLn}), (\ref{E_R2Ln-2}), and
(\ref{HL}), respectively.

Note that substituting $n=2$ either to (\ref{B_LRn-1}) or (\ref{B_RLn-1}) the
boundary $B_{LR}$ of the divergence region $D_{1/2}$ can be obtained:
\[
B_{LR}:\ \ 1-\tau_{R}\tau_{L}+\delta_{R}+\delta_{L}+\delta_{R}\delta
_{L}=0,\text{ \ }S_{0}^{LR}\in D_{L},\ S_{1}^{LR}\in D_{R}%
\]
(see Fig.~\ref{2D}), where $S_{0}^{LR}=\{(x,y):y=-\frac{\tau_{L}\delta_{R}%
}{1+\delta_{R}}x,\ x<-1\}$ and $S_{1}^{LR}=\{(x,y):y=-\frac{\tau_{R}\delta
_{L}}{1+\delta_{L}}x,\ x>-\frac{\tau_{L}}{1+\delta_{R}}\}.$

In the next section, the dynamics near regions $D_{1/5}^{R}$ and $D_{1/5}%
^{L},$ associated with WQAs, is described in more detail.

\section{Appearance of weird quasiperiodic attractors near divergence regions}

\subsection{Dynamics of map $F$ near the divergence region $D_{1/5}^{R}$}

\begin{figure} 
\includegraphics[width=1\textwidth]{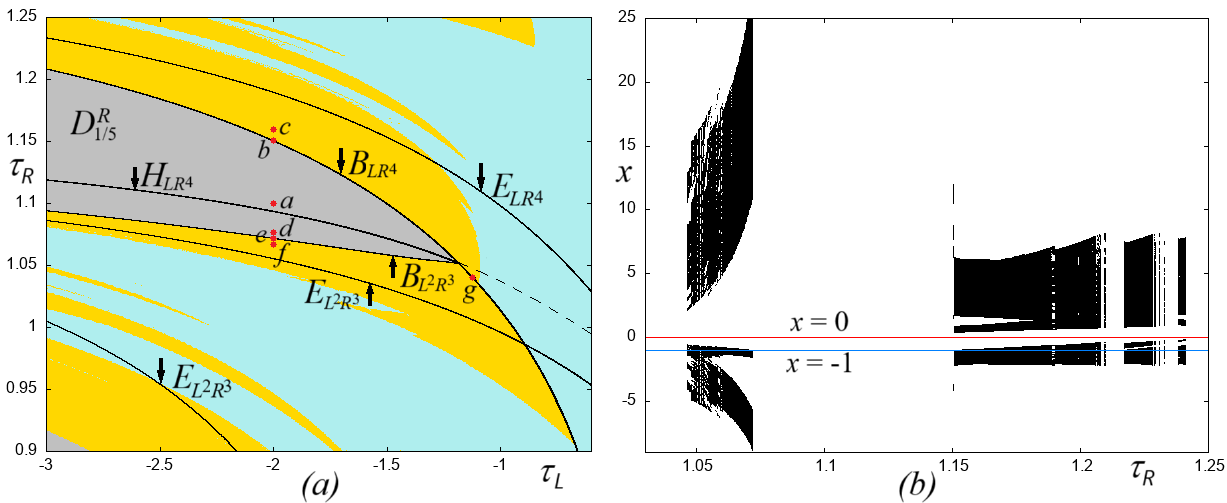}%
\caption{\label{D5R} (a) The existence regions (in yellow) of the WQAs near the divergence
region $D_{1/5}^{R}$ (in gray) bounded by curves $B_{LR^{4}}$ and
$B_{L^{2}R^{3}}$ in the $(\tau_{L},\tau_{R})$-parameter plane for 
$\delta_{L}=0.9,$ $\delta_{R}=0.7;$ (b) 1D bifurcation diagram $x$ versus $\tau_{R}$
for $\tau_{L}=-2$. The phase portrait of map $F$ at parameter points marked in
(a) by $a,$ $b,$ ..., $f$ is shown in Fig.~\ref{D5Rexam}(a),(b),...,(f),
respectively, and by $g$ in Fig.~\ref{D5Rex3}(a).} 
\end{figure}

Let us describe the dynamics of map $F$ near the divergence region
$D_{1/5}^{R}$ by considering this region in the $(\tau_{L},\tau_{R}%
)$-parameter plane for $\delta_{L}=0.9,$ $\delta_{R}=0.7$ to see how its
boundaries, $B_{LR^{4}}$ given in (\ref{B_LRn-1}) and $B_{L^{2}R^{3}}$ given
in (\ref{B_L2Rn-2}), are related to the appearance of WQAs. In Fig.
\ref{D5R}(a), the bifurcation structure near the divergence region
$D_{1/5}^{R}$ is shown magnified. Besides curves $B_{LR^{4}}$ and
$B_{L^{2}R^{3}},$ also curves $H_{LR^{4}}$, $E_{LR^{4}}$ and $E_{L^{2}R^{2}}$
given in (\ref{HR}), (\ref{E_LRn-1}), and (\ref{E_L2Rn-2}), respectively, are
shown. Note that in the considered case, the complete region $D_{1/5}^{R}$ is
inside the stability domain of the fixed point $O,$ so that in all the
examples below the fixed point $O$ is attracting. Note also that map $F$ has
$Z_{1}-Z_{0}-Z_{1}$ type of invertibility (see Property 2(a)), and it is not
difficult to show that any invariant set of map $F$ cannot have points in zone
$Z_{0}$ (i.e., in the strip between the critical lines $C^{R}$ and $C^{L}$),
or in any image of $Z_{0}$.

\begin{figure}
\includegraphics[width=1\textwidth]{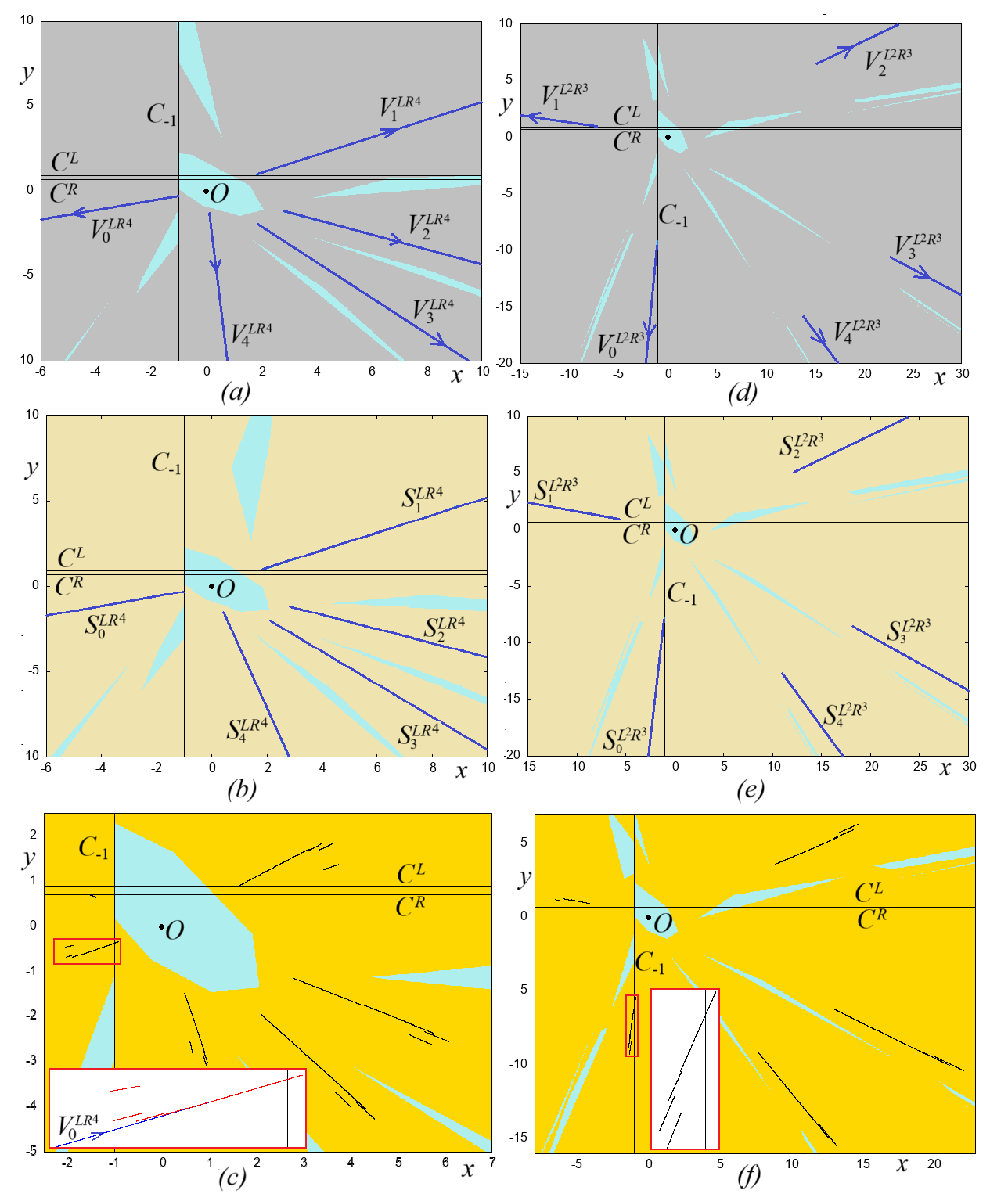} 
\caption{\label{D5Rexam} Phase portrait of map $F$ at $\delta_{L}=0.9,$ $\delta_{R}=0.7,$
$\tau_{L}=-2$ and (a) $\tau_{R}=1.1;$ (b) $\tau_{R}=1.15045;$ (c) $\tau
_{R}=1.16;$ (d) $\tau_{R}=1.075;$ (e) $\tau_{R}=1.0719;$ (f) $\tau_{R}=1.067.$
The related parameter points are marked in Fig.~\ref{D5R} by $a$, $b$, ...,
$f$, respectively. The insets in (c) and (f) show the marked windows
magnified.} 
\end{figure}

\begin{figure} 
\includegraphics[width=1\textwidth]{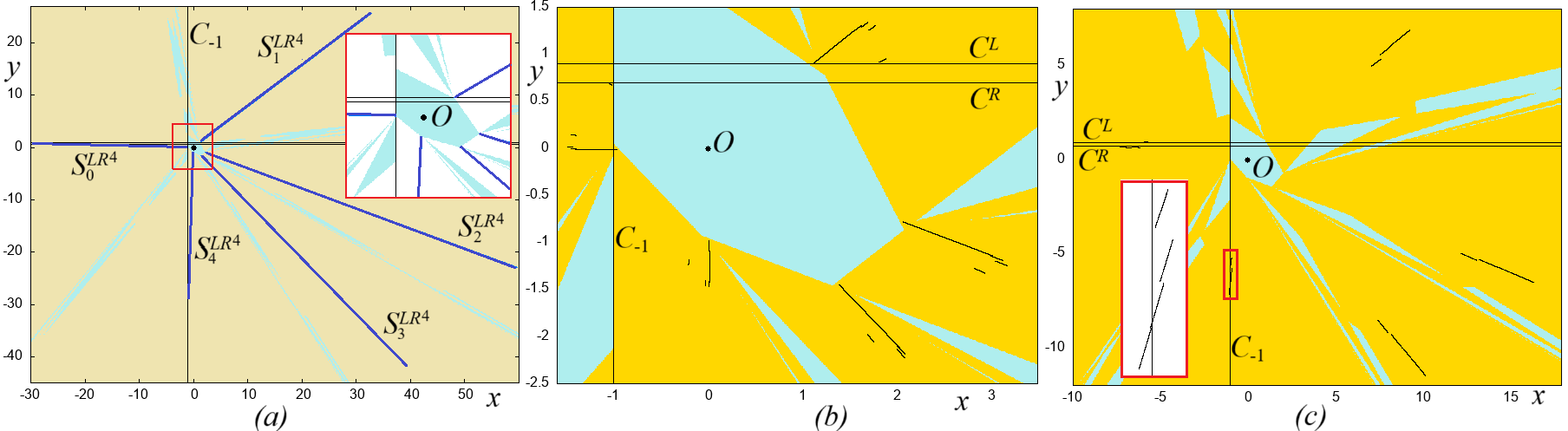} 
\caption{\label{D5Rex3} Phase portrait of map $F$ at $\delta_{L}=0.9,$ $\delta_{R}=0.7,$
$\tau_{L}=-1.125$ and (a) $\tau_{R}=1.04053,$ (b) $\tau_{R}=1.05$, (c)
$\tau_{R}=1.02$. The parameter point related to (a) is marked in Fig.
\ref{D5R}(a) by $g$. In (c), the inset shows the indicated window magnified.} 
\end{figure}

Figure \ref{D5R}(b) illustrates the dynamics by means of a 1D bifurcation
diagram showing $x$ versus $\tau_{R}$ for fixed $\tau_{L}=-2.$ At first
glance, this bifurcation diagram can be misleading because the attractors here
seem chaotic. However, as already mentioned, map $F$ cannot have chaotic
attractors, and the limit sets we observe here are WQAs born when either
boundary $B_{LR^{4}}$ (for increasing $\tau_{R}$ from a point inside the
region) or boundary $B_{L^{2}R^{3}}$ (for decreasing $\tau_{R}$) is crossed.
The mechanism of appearance of such attractors is illustrated in Fig.
\ref{D5Rexam}(a),(b),...,(f), presenting phase portraits of map $F$ at the
parameter values marked in Fig.~\ref{D5R}(a) by $a,$ $b$, ..., $f,$
respectively. One more example is shown in Fig.~\ref{D5Rex3}, associated with
the parameter point marked in Fig.~\ref{D5R}(a) by $g$. In all these figures,
the basin of the fixed point $O$ is shown in light blue, the basin of a WQA in
dark yellow, the basin of a set of segments filled with nonhyperbolic cycles
in light yellow, and the basin of diverging trajectories in gray. 

$\bullet$ Let us begin with the \textit{parameter point }$a.$ The
corresponding phase portrait of map $F$ is shown in Fig.~\ref{D5Rexam}(a). The
parameter point $a$ is located below curve $B_{LR^{4}}$ but above curve
$H_{LR^{4}}$. In this subregion of $D_{1/5}^{R},$ the attracting fixed point
$O$ coexists with divergent trajectories. More precisely, it holds that
$P_{LR^{4}}(1)<0,$ i.e., one eigenvalue of matrix $J_{LR^{4}}$ is larger than
$1$ (the second eigenvalue is obviously positive and less than $1$); moreover,
the corresponding invariant unstable halflines of map $F$ (shown in dark blue
in Fig.~\ref{D5Rexam}(a)) are located in the proper partitions: 
$V_{0}^{LR^{4}}\subset D_{L}$, $V_{j}^{LR^{4}}\subset D_{R},$ $j=\overline{1,4}.$

The left column in Fig.~\ref{D5Rexam} illustrates how the phase portrait of
map $F$ changes for increasing values of $\tau_{R},$ when the parameter point
approaches boundary $B_{LR^{4}}$ (at which $P_{LR^{4}}(1)=0$) and crosses it
(so that $P_{LR^{4}}(1)>0$ holds).

$\bullet$ \textit{Parameter point} $b.$ Fig.~\ref{D5Rexam}(b) shows the
corresponding phase portrait of map $F$, where besides the attracting fixed
point $O$ there exists a set $S^{LR^{4}}$ consisting of five cyclic halflines,
$S_{0}^{LR^{4}}\subset D_{L}$ and $S_{j}^{LR^{4}}\subset D_{R},$
$j=\overline{1,4},$ filled with nonhyperbolic $LR^{4}$-cycles, 
$S_{j}^{LR^{4}}=V_{j}^{LR^{4}}$ (see Property 6).

$\bullet$ \textit{Parameter point} $c.$ Fig.~\ref{D5Rexam}(c) presents the
phase portrait after crossing boundary $B_{LR^{4}}$, where the attracting
fixed point $O$ coexists with a WQA denoted $\mathcal{A}$, which appears after
crossing $B_{LR^{4}}$. The mechanism of creation of attractor $\mathcal{A}$
can be outlined as follows. First, note that the parameter point $c$ is above
curve $B_{LR^{4}}$ but below curve $E_{LR^{4}}$. Hence, the eigenvalues of
matrix $J_{LR^{4}}$ are both real, positive, and less than $1.$ Eigenvectors
$V_{j}^{LR^{4}},$ $j=\overline{0,4},$ which were forward invariant and
repelling before the crossing $B_{LR^{4}}$, become attracting after the
crossing. Consider halfline $V_{0}^{LR^{4}}\subset D_{L}$ (shown in blue in
the inset in Fig.~\ref{D5Rexam}(c)). Its image 
$V_{5}^{LR^{4}}=F_{LR^{4}}(V_{0}^{LR^{4}})$ necessarily intersects the discontinuity line. 
Then, under further iterations, a part of halfline $V_{5}^{LR^{4}}$ located in $D_{L},$
say, halfline $V_{5}^{LR^{4}}|_{D_{L}},$ is mapped as before by $F_{L}$ into a
halfline issuing from $C^{L}$ in $D_{R},$ while a segment of halfline
$V_{5}^{LR^{4}}$ located in $D_{R},$ say segment $V_{5}^{LR^{4}}|_{D_{R}},$ is
mapped by $F_{R}$ into a new segment issuing from $C^{R},$ which in this
example belongs to $D_{L}.$ Next iterations of segment 
$V_{5}^{LR^{4}}|_{D_{R}}$ lead to further new segments approaching $V_{j}^{LR^{4}}$,
$j=\overline{0,4}$, as the number of iterations tends to infinity (but they
are never mapped exactly into these halflines). The resulting attractor
$\mathcal{A}$ is a WQA, and it can be seen as a limit set of this iteration
process. Note that it could happen that immediately after crossing the curve
$B_{LR^{4}},$ segment $V_{5}^{LR^{4}}|_{D_{R}}$ belongs to the basin of the
attracting fixed point $O,$ and in this case a WQA is not created.

$\bullet$ Let us now turn to the \textit{parameter poin}t $d$ located above
curve $B_{L^{2}R^{3}}$ but below curve $H_{LR^{4}}.$ The related phase
portrait is shown in Fig.~\ref{D5Rexam}(d). In this case, the fixed point $O$
coexists with divergent trajectories; more precisely, it holds that
$P_{L^{2}R^{3}}(1)<0$, that is, one eigenvalue of matrix $J_{L^{2}R^{3}}$ is
larger than $1$, and the corresponding forward invariant repelling
eigenvectors of map $F$ are located in the proper partitions: $V_{0,1}%
^{L^{2}R^{3}}\subset D_{L},$ $V_{j}^{L^{2}R^{3}}\subset D_{R},$ 
$j=\overline {2,4}$ (in Fig.~\ref{D5Rexam}(d), these eigenvectors are shown in dark blue).

The right column in Fig.~\ref{D5Rexam} illustrates how the phase portrait of
map $F$ changes for decreasing values of $\tau_{R}$ when the parameter point
crosses boundary $B_{L^{2}R^{3}}.$

$\bullet$ \textit{Parameter point} $e.$ Fig.~\ref{D5Rexam}(e) shows the
related phase portrait. It holds that $(\tau_{L},\tau_{R})\in B_{L^{2}R^{3}}$,
thus $P_{L^{2}R^{3}}(1)=0,$ and besides the attracting fixed point $O,$ there
exists a set $S^{L^{3}R^{3}}$ consisting of $5$-cyclic halflines,
$S_{0,1}^{L^{3}R^{3}}\subset D_{L}$ and $S_{j}^{L^{2}R^{3}}\subset D_{R},$
$j=\overline{2,4}$ ($S^{L^{3}R^{3}}=V_{j}^{L^{2}R^{3}}$), filled with
nonhyperbolic $L^{2}R^{3}$-cycles (see Property 6);

$\bullet$ \textit{Parameter point} $f.$ Fig.~\ref{D5Rexam}(f) presents the
corresponding phase portrait. The parameter point $f$ is below 
$B_{L^{2}R^{3}}$, so that $P_{L^{2}R^{3}}(1)>0,$ and the attracting fixed point $O$ coexists
with a WQA denoted $\mathcal{A}$,\ born after crossing boundary 
$B_{L^{2}R^{3}}$. The mechanism of its creation is similar to the one described in the
previous example. Note that in this case, the parameter point is below curve
$B_{L^{2}R^{3}}$ but above curve $E_{L^{2}R^{3}}$. Hence, the eigenvalues of
matrix $J_{L^{2}R^{3}}$ are both real, positive, and less than $1$.

It is interesting to compare the WQAs shown in Figs. \ref{D5Rexam}(c) and
\ref{D5Rexam}(f). Considering map $F^{5},$ it can be determined that the
attractor in Fig.~\ref{D5Rexam}(f) consists of five \textit{cyclic} blocks,
unlike the attractor shown in Fig.~\ref{D5Rexam}(c), which is not cyclic. In
fact, the complete block shown magnified in the inset in Fig.~\ref{D5Rexam}(f)
is mapped in one block in partition $D_{L},$ while the block of the attractor
shown magnified in the inset in Fig.~\ref{D5Rexam}(c) is mapped partially in
$D_{L}$ and partially in $D_{R},$ breaking the cyclicity of this attractor.

$\bullet$ \textit{Parameter point} $g$. This example is related to the case
when the $(\tau_{L},\tau_{R})$-parameter point belongs to curve $B_{LR^{4}}.$
Since this point is below curve $H_{LR^{4}},$ the $LR^{4}$-cycle at infinity
does not exist and curve $B_{LR^{4}}$ no longer serves as a boundary of the
divergence region $D_{1/5}^{R}.$ However, an admissible set $S^{LR^{4}}$ still
exists: it consists of five cyclic \textit{bounded} segments, 
$S_{0}^{LR^{4}}\subset D_{L}$ and $S_{j}^{LR^{4}}\subset D_{R},$ $j=\overline{1,4},$ filled
with nonhyperbolic $LR^{4}$-cycles. Figure \ref{D5Rex3}(a) presents the phase
portrait of map $F$, where the attracting fixed point $O$ coexists with
segments $S_{j}^{LR^{4}},$ $j=\overline{0,4},$ shown in dark blue.

Starting from the case show in Fig.~\ref{D5Rex3}(a), a WQA denoted
$\mathcal{A}$ appears for both increasing and decreasing values of $\tau_{R}$,
as in Fig.~\ref{D5Rex3}(b) and Fig.~\ref{D5Rex3}(c), respectively. In
particular, after crossing $B_{LR^{4}}$ for increasing $\tau_{R}$ (see Fig.
\ref{D5Rex3}(b)), eigenvectors $V_{j}^{LR^{4}}$ become attracting, so that
segment $F^{5}(V_{0}^{LR^{4}})$ necessarily intersects the discontinuity line
$C_{-1}$ and the mechanism of creation of attractor $\mathcal{A}$ is similar
to the one illustrated in Figs. \ref{D5Rexam}(b) and (c). Note that for
further increasing $\tau_{R},$ attractor $\mathcal{A}$ disappears quite
quickly due to a contact with the basin boundary of the fixed point $O$. In
contrast, crossing $B_{LR^{4}}$ for decreasing $\tau_{R}$ (see Fig.
\ref{D5Rex3}(c)), eigenvectors $V_{j}^{LR^{4}}$ become repelling, so that
segment $F^{5}(V_{4}^{LR^{4}})$ necessarily intersects the discontinuity line
$C_{-1}$, and the mechanism of creation of a WQA is associated with further
images of segment $F^{5}(V_{4}^{LR^{4}})|_{D_{L}}.$%

\subsection{Dynamics of map $F$ near the divergence region $D_{1/5}^{L}$}

\begin{figure}
\includegraphics[width=1\textwidth]{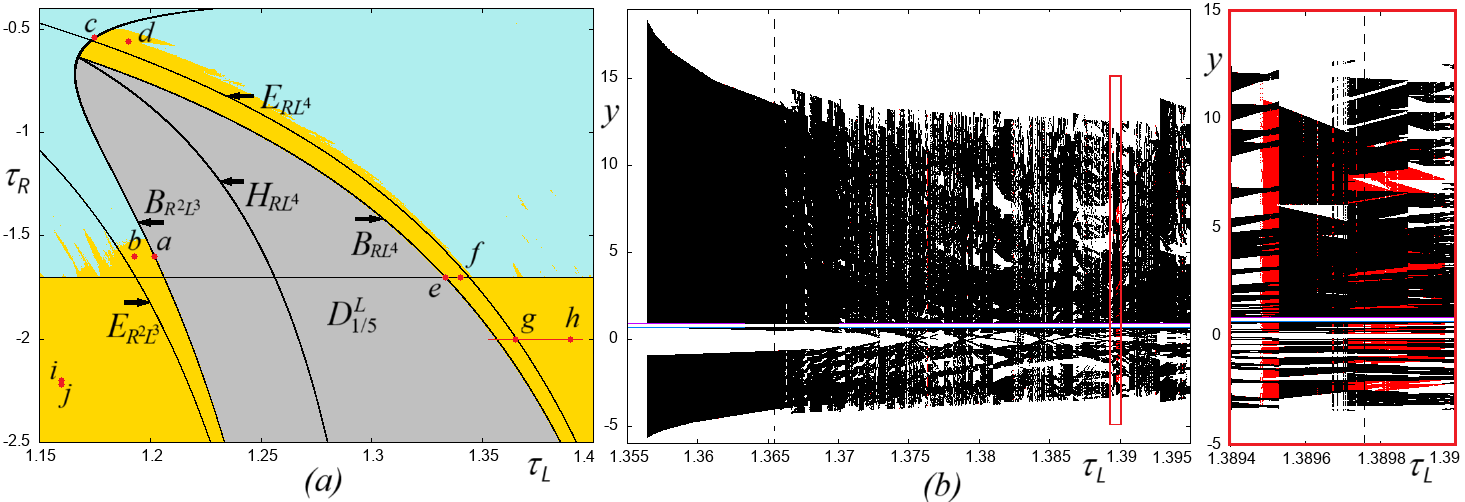} 
\caption{\label{D5L} (a) The existence regions (in yellow) of the WQAs near the divergence
region $D_{1/5}^{L}$ (in gray) bounded by curves $B_{RL^{4}}$ and
$B_{R^{2}L^{3}}$ in the $(\tau_{L},\tau_{R})$-parameter plane for 
$\delta_{L}=0.9,$ $\delta_{R}=0.7;$ (b) 1D bifurcation diagram $y$ versus $\tau_{L}$
for $\tau_{R}=-2$. The rightmost figure presents a magnified window indicated
in (b), where coexisting WQAs (shown in red and black) are better visible. The
phase portraits of map $F$ at parameter points marked by $a,$ $b,$ ..., $j$
are shown in Figs. \ref{D5Lex1}(a),(b), \ref{D5Lex2}(a),(b), \ref{D5Lex3}%
(a),(b), \ref{D5Lex4}(a),(b), and \ref{D5Lex5}(a),(b), respectively.} 
\end{figure}

We now turn to an example of the bifurcation structure near the divergence
region $D_{1/5}^{L}$. Figure \ref{D5L}(a) shows the $(\tau_{L},\tau_{R}%
)$-parameter plane for $\delta_{L}=0.9,$ $\delta_{R}=0.7$, where besides
curves $B_{RL^{4}}$ and $B_{R^{2}L^{3}},$ given in (\ref{B_RLn-1}) and
(\ref{B_R2Ln-2}), also curves $H_{RL^{4}}$, $E_{RL^{4}},$ and 
$E_{R^{2}L^{3}},$ given in (\ref{HL}), (\ref{E_RLn}), and (\ref{E_R2Ln-2}), 
respectively,
are shown. Comparing with the divergence region $D_{1/5}^{R}$ in Fig.
\ref{D5R}(a), region $D_{1/5}^{L}$ in Fig.~\ref{D5L}(a) is located not
entirely in the stability domain of the fixed point $O,$ but also in the
region below the straight line $\tau_{R}=-1.7,$ where the fixed point $O$ is a
saddle. Figure \ref{D5L}(b) and its magnified part in the right panel present
a 1D bifurcation diagram for $y$ versus $\tau_{L},$ setting $\tau_{R}=-2.$ The
corresponding parameter path is indicated in Fig.~\ref{D5L}(a) by the red
horizontal segment. This diagram confirms, in particular, that attractors of
map $F$ have no points in the strip between the critical lines $C^{R}$ and
$C^{L}$ (zone $Z_{0}$ in Property 2(a)), defined in the considered case by
$y=0.7$ and $y=0.9$, respectively (see two horizontal straight lines in Fig.
\ref{D5L}(b)). Below, we give a few examples related to Fig.~\ref{D5L}(b).

\begin{figure}
\includegraphics[width=1\textwidth]{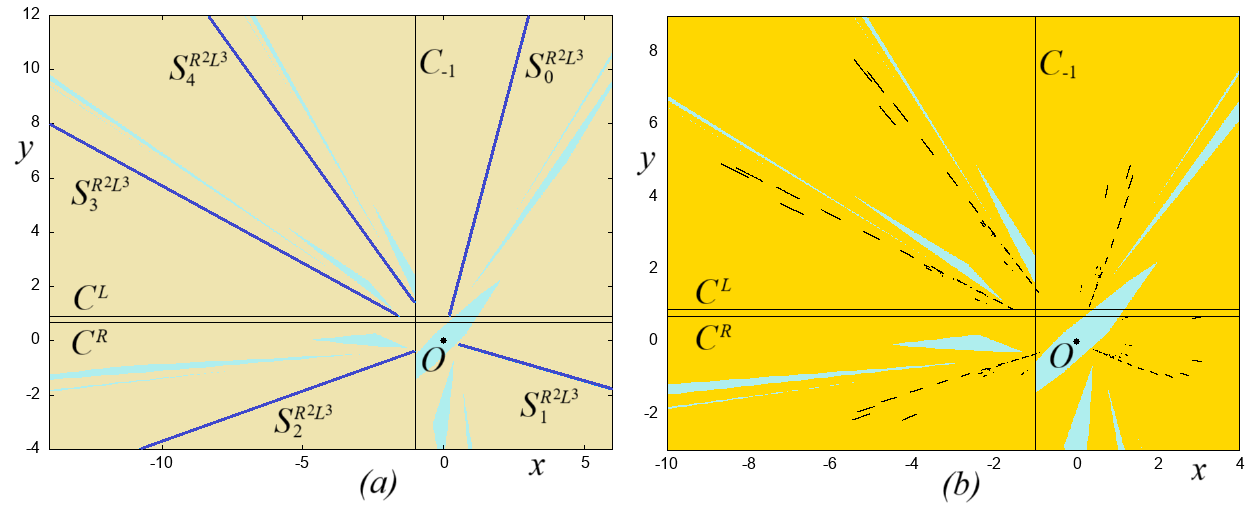}
\caption{\label{D5Lex1} Phase portrait of map $F$ at $\delta_{L}=0.9,$ $\delta_{R}=0.7,$
$\tau_{R}=-1.6$ and (a) $\tau_{L}=1.201945,$ (b) $\tau_{L}=1.193.$ 
In Fig.~\ref{D5L}(a), the related parameter points are marked by $a$ and $b$,
respectively.} 
\end{figure}

To illustrate the dynamics of map $F$ near the divergence region 
$D_{1/5}^{L},$ we present the phase portraits of map $F$ for parameter values
indicated in Fig.~\ref{D5L}(a) by $a,$ $b,$ $c,$...$,$ $j.$ The following
should be noted:

$\bullet$ \textit{Parameter points} $a$ \textit{and} $b$. In Fig.
\ref{D5Lex1}(a), related to the parameter point $a$ belonging to curve
$B_{R^{2}L^{3}}$, besides the attracting fixed point $O,$ there exists a set
$S^{R^{2}L^{3}}$ consisting of five cyclic halflines, 
$S_{0,1}^{R^{2}L^{3}}\subset D_{R}$ and $S_{j}^{R^{2}L^{3}}\subset D_{L},$ 
$j=\overline{2,4}$ ($S_{j}^{R^{2}L^{3}}=V_{j}^{R^{2}L^{3}}$), filled with nonhyperbolic
$R^{2}L^{3}$-cycles (see Property 6). Unlike crossing boundary $B_{RL^{4}}$,
crossing boundary $B_{R^{2}L^{3}}$ may not lead to a WQA. In fact, after
crossing $B_{RL^{4}},$ eigenvectors $V_{j}^{R^{2}L^{3}},$ $j=\overline{0,4},$
become attracting, thus after five iterations, halfline 
$F^{5}(V_{4}^{R^{2}L^{3}})$ and/or $F^{5}(V_{2}^{R^{2}L^{3}})$ necessarily intersects
partition $D_{R}.$ It may happen that the segment(s) of 
$F^{5}(V_{4}^{R^{2}L^{3}})$ and/or $F^{5}(V_{2}^{R^{2}L^{3}})$ located in $D_{R}$
immediately after the bifurcation enter(s) the basin of the attracting fixed
point $O,$ and thus no WQA appears. It can be seen in Fig.~\ref{D5Lex1}(a)
that segments $S_{j}^{R^{2}L^{3}}$ are relatively far from the basin of the
fixed point $O,$ and for decreasing $\tau_{L}$ a WQA, as, e.g., in Fig.
\ref{D5Lex1}(b), is created. Note that the related parameter point $b$ in Fig.
\ref{D5L}(a) is located between curves $B_{R^{2}L^{3}}$ and $E_{R^{2}L^{3}},$
where the eigenvalues of the matrix $J_{R^{2}L^{3}}$ are both real, positive,
and less than $1$.

\begin{figure} 
\includegraphics[width=1\textwidth]{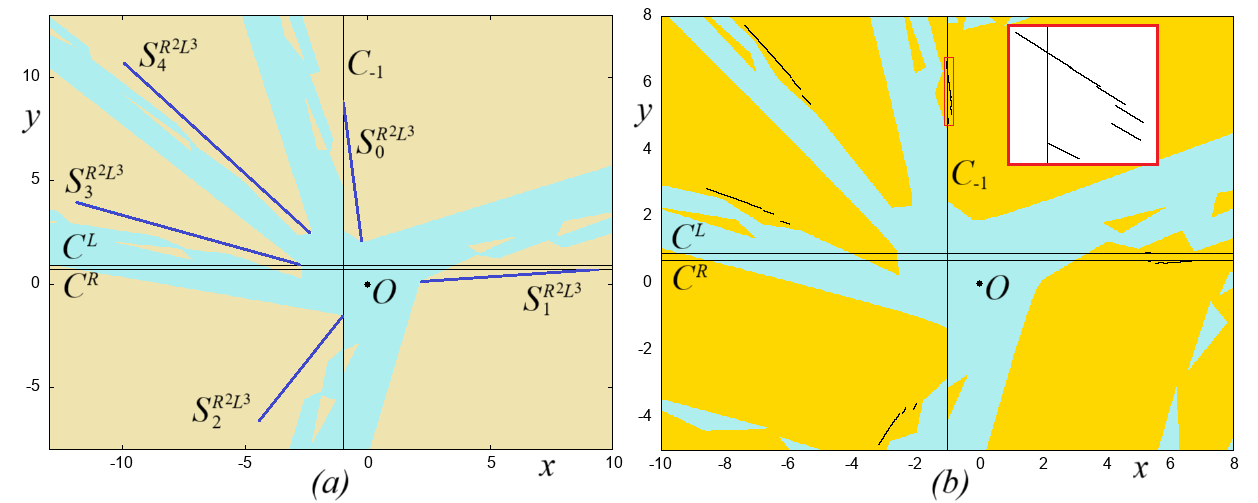}
\caption{\label{D5Lex2} Phase portrait of map $F$ at $\delta_{L}=0.9,$ $\delta_{R}=0.7$ and
(a) $\tau_{L}=1.175,$ $\tau_{R}=-0.543,$ (b) $\tau_{L}=1.19,$ $\tau
_{R}=-0.56.$ In Fig.~\ref{D5L}(a), the related parameter points are marked by
$c$ and $d$, respectively.} 
\end{figure}

$\bullet$ \textit{Parameter points} $c$ \textit{and} $d$. In Fig.
\ref{D5Lex2}(a), the phase portrait of map $F$ is shown for the parameter
point $c$ belonging to curve $B_{R^{2}L^{3}}$. In contrast to the case in Fig.
\ref{D5Lex1}(a), the parameter point is above curve $H_{RL^{4}}.$ In this
region, an $R^{2}L^{3}$-cycle at infinity does not exist, and this part of
curve $B_{R^{2}L^{3}}$ is not at the boundary of the divergence region
$D_{1/5}^{L}.$ However, an admissible set $S^{R^{2}L^{3}}$ still exists,
consisting of five cyclic \textit{bounded} segments, 
$S_{0,1}^{R^{2}L^{3}}\subset D_{R}$ and $S_{j}^{R^{2}L^{3}}\subset D_{L},$ $j=\overline{2,4}$
($S_{j}^{R^{2}L^{3}}=V_{j}^{R^{2}L^{3}}$) filled with nonhyperbolic
$R^{2}L^{3}$-cycles (see Property 6). An example of a WQA existing for the
parameter values on the right of curve $B_{R^{2}L^{3}}$ is shown in Fig.
\ref{D5Lex2}(b). Note that the related parameter point $d$ in Fig.
\ref{D5L}(a) belongs to the region where eigenvectors $V_{j}^{R^{2}L^{3}},$
$j=\overline{0,4},$ are repelling. Comparing with the example in Fig.
\ref{D5Lex1}, boundary $B_{R^{2}L^{3}}$ is crossed in the opposite direction:
segment $F^{5}(V_{0}^{R^{2}L^{3}})$ necessarily intersects partition $D_{L}$
(see the inset in Fig.~\ref{D5Lex2}(b)). Further iterations of segment
$F^{5}(V_{0}^{R^{2}L^{3}})|_{D_{L}}$ lead to new segments, and the WQA can be
seen as a limit set of this process.

\begin{figure}
\includegraphics[width=1\textwidth]{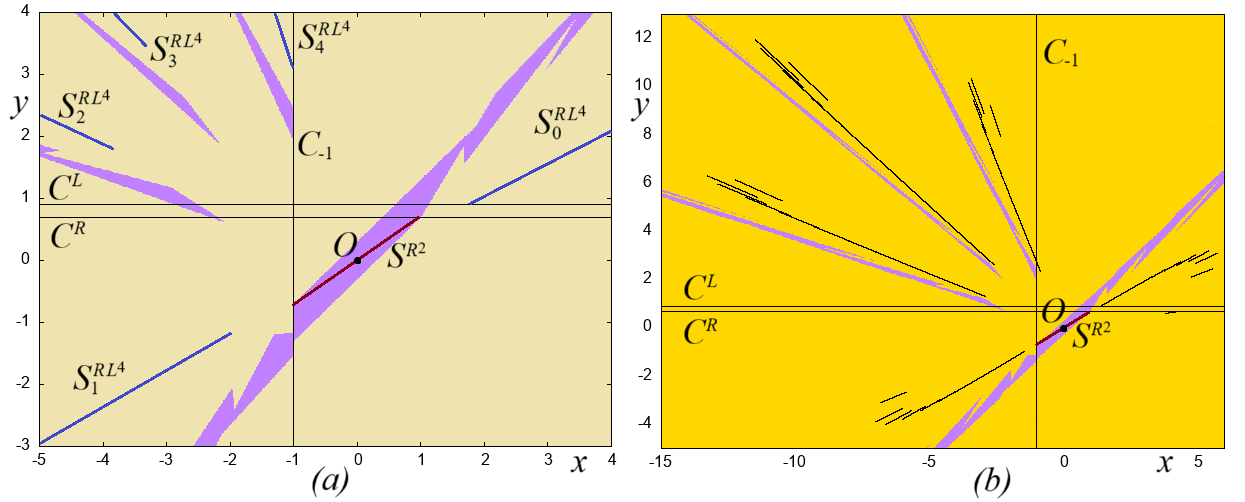}
\caption{\label{D5Lex3} Phase portrait of map $F$ at $\delta_{L}=0.9,$ $\delta_{R}=0.7,$
$\tau_{R}=-1.7$ and (a) $\tau_{L}=1.333153,$ (b) $\tau_{L}=1.34.$ In Fig.
\ref{D5L}(a), the related parameter points are marked by $e$ and $f$,
respectively.} 
\end{figure}

$\bullet$ \textit{Parameter points} $e$ \textit{and} $f.$ In Fig.
\ref{D5Lex3}(a) a special case is illustrated where parameter point $e$
belongs to curve $B_{RL^{4}}$ and $\tau_{R}=-1.7$. Recall that at 
$\tau_{R}=-1.7,$ the fixed point $O$ undergoes a degenerate flip bifurcation (see
Property 4): there is a segment $S^{R^{2}}=\{(x,y):-1<x<1,\ y=0.7x\}$ filled
with points of nonhyperbolic $2$-cycles (except the fixed point $O$), and this
segment coexists with a set $S^{RL^{4}}$ consisting of five cyclic halflines,
$S_{0}^{RL^{4}}\subset D_{R}$ and $S_{j}^{RL^{4}}\subset D_{L},$
$j=\overline{1,4},$ filled with nonhyperbolic $RL^{4}$-cycles. An example of a
WQA born on the right side of curve $B_{RL^{3}}$ and coexisting with segment
$S^{R^{2}}$ is shown in Fig.~\ref{D5Lex3}(b). The related parameter point is
marked $f$ in Fig.~\ref{D5L}(a).

\begin{figure}
\includegraphics[width=1\textwidth] {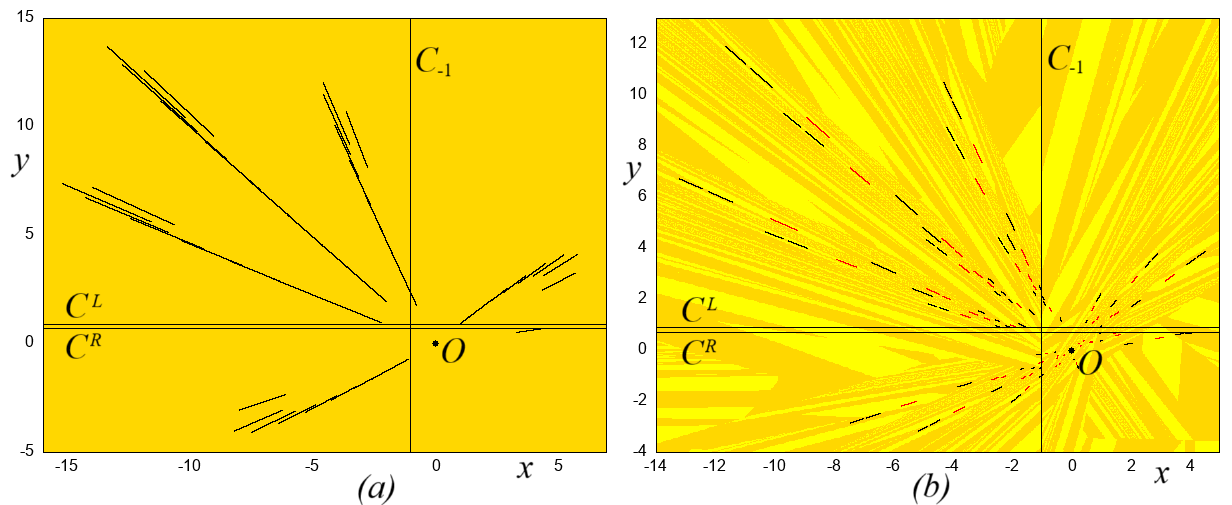}
\caption{\label{D5Lex4} Phase portrait of map $F$ at $\delta_{L}=0.9,$ $\delta_{R}=0.7,$
$\tau_{R}=-2$ and (a) $\tau_{L}=1.365,$ (b) $\tau_{L}=1.38975.$ In Fig.
\ref{D5L}(a), the related parameter points are marked by $g$ and $h$,
respectively.} 
\end{figure}

$\bullet$ \textit{Parameter points} $g$ \textit{and} $h$. In Fig.
\ref{D5Lex4}, two examples of the phase portrait of map $F$ are shown related
to the 1D bifurcation diagram in Fig.~\ref{D5L}(b). It is easy to verify that
the WQA shown in Fig.~\ref{D5Lex4}(a) related to the parameter point $g$,
consists of five cyclic blocks, that is, each block is $F^{5}$-invariant. In
this case, only one block intersects the discontinuity line, and this entire
block is mapped into the same partition $D_{R}.$ The cyclicity is not lost
immediately but soon after crossing curve $E_{RL^{4}}$ (at 
$\tau_{L}\approx1.36540484,$ see the dashed vertical line in Fig.~\ref{D5L}(b)), when
the eigenvalues of matrix $J_{RL^{4}}$ become complex-conjugate. Figure
\ref{D5Lex4}(b) related to the parameter point $h$, illustrates two coexisting
WQAs and their basins at $\tau_{L}=1.38975$ (see the dashed vertical line in
the right panel of Fig.~\ref{D5L}(b)). As can be seen in the right panel of
Fig.~\ref{D5L}(b), there are several intervals of bistability, confirming that
the coexistence of two WQAs may be persistent under parameter perturbations.

\begin{figure} 
\includegraphics[width=1\textwidth]{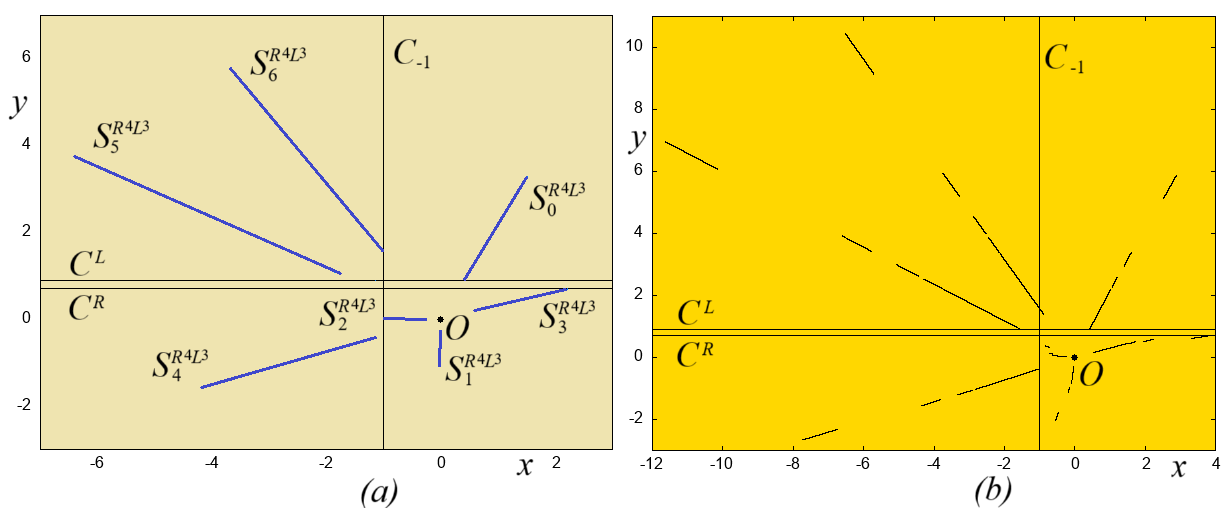}
\caption{\label{D5Lex5} Phase portrait of map $F$ at $\delta_{L}=0.9,$ $\delta_{R}=0.7,$
$\tau_{L}=1.16$ and (a) $\tau_{R}=-2.199528,$ (b) $\tau_{L}=-2.219.$ In Fig.
\ref{D5L}(a), the related parameter points are marked by $i$ and $j$,
respectively.} 
\end{figure}

$\bullet$ \textit{Parameter points} $i$ \textit{and} $j.$ Finally, it is worth
noting that in the $(\tau_{L},\tau_{R})$-parameter plane in Fig.~\ref{D5L}(a),
there are parameter values associated with segments filled with nonhyperbolic
cycles having symbolic sequences that are different from those considered
above, i.e., from the basic and complementary symbolic sequences. An example
is shown in Fig.~\ref{D5Lex5}(a), in the case when map $F$ has a set of seven
cyclic segments $S_{j}^{R^{4}L^{3}},$ $j=\overline{0,6},$ filled with
nonhyperbolic $R^{4}L^{3}$-cycles (the related point is marked $i$). A WQA
appearing for decreasing $\tau_{R},$ e.g., at the parameter point $j$, is
shown in Fig.~\ref{D5Lex5}(b). A mechanism of creation of this attractor is
associated with segment $V_{6}^{R^{4}L^{3}}\subset D_{L},$ which now, after
seven iterations, intersects the discontinuity line that leads to a new
segment in $D_{R},$ whose images approach, as $n\rightarrow\infty,$ segments
$V_{j}^{R^{4}L^{3}}.$ The WQA is a limit set of this process. This is not
surprising, since we have shown only a few curves in the parameter space
associated with $P_{\sigma}(1)=0,$ but clearly there are infinitely many such
curves (although they all have zero Lebesgue measure), and they may be
associated with the appearance of WQAs.

\section{Conclusions}

We considered a 2D piecewise linear map defined by two linear homogeneous
functions in two partitions of the phase plane separated by a vertical
discontinuity line. This map can have a new type of attractor, which we call a
weird quasiperiodic attractor. This type of attractor was first observed in
applied models (see, e.g., \cite{Kollar}, \cite{Gardini et al 25},
\cite{Gardini et al 25c}). A WQA may consist of an infinite number of
segments, fill some 2D subsets of the phase plane in a (seemingly) dense
manner, or have a mixed structure. At first glance, a WQA may be
misinterpreted as a chaotic attractor due to its intricate (weird) structure.
However, chaos cannot be observed in the considered map since it cannot have
hyperbolic cycles. The complex structure of the attractor is primarily related
to the discontinuity of the map and the homogeneity of its linear components.

The present paper can be seen as a first step towards understanding the basic
properties of WQAs and the mechanisms of their emergence. We described the
bifurcation structure of the parameter space of map $F$ and analytically
obtained some boundaries of the divergence regions, crossing which a WQA can
appear. Special attention was paid to the case when map $F$ has $Z_{1}%
-Z_{0}-Z_{1}$ type of invertibility (an analog of a 1D gap map). To illustrate
how a WQA can appear/disappear, we described in detail the corresponding
transformations of the phase portrait of the map.

However, many issues remain unresolved. The first set of questions concerns
the class of maps that can have WQAs. It is clear that this class includes
discontinuous piecewise linear maps of dimension two or higher, defined by two
or more homogeneous linear functions, or linear functions with the same fixed
point (which can be translated to the origin, leading to the homogeneous
case). It is also clear that discontinuity is a necessary element for the
existence of a WQA. However, it remains an open problem whether linearity of
the map components is a necessary condition or whether the functions defining
the map can be nonlinear. There are also questions related to the structure of
WQAs. For example, we conjecture that the structure of a WQA cannot be
fractal, but it is not easy to prove this statement in the general case. One
of the main questions is related to the quasiperiodicity of the trajectory on
the attractor. It is well known that for the standard quasiperiodic attractor,
such as an attracting closed invariant curve in the case of irrational
rotation, the maximum Lyapunov exponent is zero, and there is no sensitive
dependence on initial conditions. In the case of WQA, dependence on initial
conditions is inevitable due to the discontinuity of the map. Taking two
adjacent points belonging to the attractor but on opposite sides of the
discontinuity line will necessarily lead to two images of these points that
are distant from each other. In \cite{Kollar}, this property is called weak
sensitive dependence on initial conditions. The question is how weak this
dependence is and whether the maximum Lyapunov exponent remains zero
(numerically it is slightly positive, but we think this may be a numerical
problem). It is also worth mentioning an interesting problem related to the
mechanism of the emergence of WQAs and their properties when the map is
noninvertible of $Z_{1}-Z_{2}-Z_{1}$ type. Some of the aforementioned problems
are addressed in our companion paper \cite{Gardini et al 25b}, while others
are left for further research.\bigskip

\textbf{Acknowledgements}

Laura Gardini thanks the Czech Science Foundation (Project 22-28882S), the
VSB---Technical University of Ostrava (SGS Research Project SP2024/047), the
European Union (REFRESH Project-Research Excellence for Region Sustainability
and High-Tech Industries of the European Just Transition Fund, Grant
CZ.10.03.01/00/22 003/000004). Davide Radi thanks the Gruppo Nazionale di
Fisica Matematica GNFM-INdAM for their financial support. The work of Davide
Radi and Iryna Sushko was funded by the European Union - Next Generation EU,
Mission 4: "Education and Research" -- Component 2: "From research to
business", through PRIN 2022 under the Italian Ministry of University and
Research (MUR). Project: 2022JRY7EF -- Qnt4Green -- Quantitative Approaches
for Green Bond Market: Risk Assessment, Agency Problems and Policy Incentives
- CUP: J53D23004700008.

\end{document}